\documentclass[a4paper,12pt]{article}
\usepackage{amsmath, amssymb, latexsym, amscd, amsthm,amsfonts,amstext}
\usepackage[mathscr]{eucal}
\usepackage{graphicx}
\usepackage{subfig}

\numberwithin{equation}{section}
\begin{document}

\title{ Spatiotemporal Monitoring of Epidemics via Solution of a Coefficient Inverse Problem
\thanks{The third author was supported in part by China Postdoctoral Science Foundation: 2023M731528.}
}

\author{ Michael V. Klibanov\footnotemark[2]  \and Jingzhi Li \footnotemark[3] \and Zhipeng Yang\footnotemark[4]}
\date{}

\maketitle

\renewcommand{\thefootnote}{\fnsymbol{footnote}}

\footnotetext[2]{Department of Mathematics and Statistics, University of North \and Carolina at Charlotte, Charlotte, NC 28223, USA, (mklibanv@uncc.edu).}

\footnotetext[3]{Department of Mathematics \& National Center for Applied
	Mathematics Shenzhen \& SUSTech International Center for Mathematics,
	Southern University of Science and Technology, Shenzhen 518055, P.~R.~China, (li.jz@sustech.edu.cn).}

\footnotetext[4]{Department of Mathematics, Southern University of Science and Technology, Shenzhen 518055, P. R. China, (yangzp@sustech.edu.cn).}

\renewcommand{\thefootnote}{\arabic{footnote}}

\date{}
\maketitle

\begin{abstract}
Let S,I and R be susceptible, infected and recovered populations in a city
affected by an epidemic. The SIR model of Lee, Liu, Tembine, Li and Osher,
\emph{SIAM J. Appl. Math.},~81, 190--207, 2021 of the spatiotemoral spread
of epidemics is considered. This model consists of a system of three
nonlinear coupled parabolic Partial Differential Equations with respect to
the space and time dependent functions S,I and R. For the first time, a
Coefficient Inverse Problem (CIP) for this system is posed. The so-called
\textquotedblleft convexification" numerical method for this inverse problem
is constructed. The presence of the Carleman Weight Function (CWF) in the
resulting regularization functional ensures the global convergence of the
gradient descent method of the minimization of this functional to the true
solution of the CIP, as long as the noise level tends to zero. The CWF is
the function, which is used as the weight in the Carleman estimate for the
corresponding Partial Differential Operator. Numerical studies demonstrate
an accurate reconstruction of unknown coefficients as well as S,I,R
functions inside of that city. As a by-product, uniqueness theorem for this
CIP is proven. Since the minimal measured input data are required, then the
proposed methodology has a potential of a significant decrease of the cost
of monitoring of epidemics.
\end{abstract}



\section{Introduction}

\label{sec:1}

The experience of the COVD-19 pandemic demonstrates, in particular, the
importance of the mathematical modeling of the spread of epidemics. The
commonly used model is the model of Kermack and McKendrick, which was
proposed in 1927 \cite{Kermack}. However, since this model is based on a
system of three coupled Ordinary Differential Equations, then it describes
only the total number of susceptible (S), infected (I) and recovered (R)
populations, the so-called \textquotedblleft SIR model". On the other hand,
recently Lee, Liu, Tembine, Li and Osher \cite{Lee} have extended the SIR
model of \cite{Kermack} to the case of a system of three nonlinear coupled
parabolic Partial Differential Equations (PDEs), which we call again
\textquotedblleft SIR model". The advantage of \cite{Lee} over \cite{Kermack}
is that the model of \cite{Lee} governs both spatial and time dependencies
of SIR populations. Thus, it makes sense to use this model for monitoring of
both pointwise and timewise distributions of SIR populations. To decrease
the cost, it is desirable to decrease the number of measurements. Results of
the current paper indicate that it is possible to achieve this goal.

Let $\mathbf{x}=\left( x,y\right) \in \mathbb{R}^{2}$ be the vector of
spatial variables and $t>0$ be time. The system of PDEs of \cite{Lee}
contains some parameters, which are actually unknown functions. These
parameters are the infection rate $\beta \left( \mathbf{x},t\right) ,$ the
recovery rate $\gamma \left( \mathbf{x},t\right) $ and three 2D vector
functions $q_{S}\left( \mathbf{x},t\right) ,q_{I}\left( \mathbf{x},t\right) $
and $q_{R}\left( \mathbf{x},t\right) ,$ which are velocities of propagations
of S, I and R populations respectively. Therefore, these functions are
subjects to the solutions of certain Coefficient Inverse Problems (CIPs).
This is the first work concerning a CIP for the SIR system of \cite{Lee}.
Hence, as in any first work about an applied problem, it is natural to
introduce some simplifications. More precisely, we assume that the unknown
coefficients $\beta $ and $\gamma $ depend only on $\mathbf{x}$, i.e. $\beta
=\beta \left( \mathbf{x}\right) $ and $\gamma =\gamma \left( \mathbf{x}%
\right) .$ In addition, we assume that the 2D vector functions $q_{S}\left(
\mathbf{x}\right) ,q_{I}\left( \mathbf{x}\right) $,$q_{R}\left( \mathbf{x}%
\right) $ of velocities are known and depend only on $\mathbf{x}$. More
complicated cases can be considered later. We assume that
\begin{equation}
	\beta \left( \mathbf{x}\right) ,\gamma \left( \mathbf{x}\right) \in C\left(
	\overline{\Omega }\right) ;\text{ }q_{S}\left( \mathbf{x}\right)
	,q_{I}\left( \mathbf{x}\right) ,q_{R}\left( \mathbf{x}\right) \in
	C^{1}\left( \overline{\Omega }\right) ,  \label{1.0}
\end{equation}%
where $\Omega \subset \mathbb{R}^{2}$ is a bounded domain occupied by a city
affected by an epidemic. Below $\partial \Omega $ is the piecewise smooth
boundary of $\Omega .$ Thus, the CIP of this paper targets the simultaneous
reconstruction of all components of the 5D vector function $\Phi \left(
\mathbf{x},t\right) ,$
\begin{equation}
	\Phi \left( \mathbf{x},t\right) =\left( \beta \left( \mathbf{x}\right)
	,\gamma \left( \mathbf{x}\right) ,S\left( \mathbf{x},t\right) ,I\left(
	\mathbf{x},t\right) ,R\left( \mathbf{x},t\right) \right) .  \label{1.1}
\end{equation}

The input data are the minimal ones. This is a single measurement of the 3D
vector function $\Psi \left( \mathbf{x}\right) =\left( S\left( \mathbf{x}%
,t_{0}\right) ,I\left( \mathbf{x},t_{0}\right) ,R\left( \mathbf{x}%
,t_{0}\right) \right) $ inside of that city at a single fixed moment of time
$\left\{ t=t_{0}>0\right\} $ as well as boundary measurements at $\partial
\Omega $ of functions $S\left( \mathbf{x},t\right) ,I\left( \mathbf{x}%
,t\right) ,R\left( \mathbf{x},t\right) $ and their normal derivatives.
Normal derivatives are fluxes of SIR populations through the boundary of
that city. In particular, fluxes can be assumed to be zeros, as in \cite{Lee}%
. We do not impose this assumption in our theory, although we use it in
computations. It is hardly practical to measure the S,I,R functions at $%
\left\{ t=0\right\} .$ Indeed, the epidemic process is yet immature at small
times, and the authorities do not even know about the existence of an
epidemic at $t\approx 0.$

Currently measurements of SIR populations are carried out at many locations
inside of a city affected by an epidemic, and at all times $t\in \left(
t_{1},T\right) ,$ where $t_{1}>0$ is a certain moment of time when
authorities start to fight that epidemic, and $T>t_{1}$ is also a certain
moment of time. On the other hand, our technique provides an accurate
reconstruction of the vector function $\Phi $ in (\ref{1.1}) inside of the
city for all times of interest by requiring measurements inside of that city
only at a single moment of time $t_{0}\in \left( t_{1},T\right) $ and
conducting the rest of the measurements only at the city's boundary. Hence,
our technique has a potential of a significant decrease of the cost of
monitoring of epidemics.

All CIPs for PDEs are both ill-posed and nonlinear ones. The conventional
approach to numerical solutions of CIPs is based on the optimization of a
least squares mismatch functional, see, e.g. \cite%
{B1,B2,B3,B4,Chavent,Gonch1,Gonch2}. However, phenomena of the ill-posedness
and the nonlinearity of CIPs cause the non-convexity of those functionals.
In turn, this can lead to the presence of multiple local minima and ravines,
posing a risk of these optimization-based methods becoming stuck, see, e.g.
\cite{Scales} for a good numerical example of multiple local minima.
Consequently, to effectively compute numerical solutions with these
conventional techniques, good initial guesses are necessary. However, this
requirement is not often met in practical situations.

In this work we develop the so-called \textquotedblleft convexification"
numerical method for our CIP. The convexification is the concept for CIPs,
which is based on Carleman estimates. This concept was originally designed
in \cite{KI,Klib97} to avoid multiple local minima and ravines of
conventional least squares mismatch functionals. Each new CIP requires its
own version of the convexification. While the originating works \cite%
{KI,Klib97} where purely theoretical ones, more recently this research team
has published a number of works, in which the analytical studies of several
versions of the convexification are combined with numerical results for a
variety of CIPs, see, e.g. \cite{KL}-\cite{MFG8}. In addition, the second
generation of the convexification was developed, see, e.g. \cite%
{Baud1,Baud2,Baud3,NguyenKl2022,NguyenV}.

The convexification consists of two steps:

\begin{enumerate}
	\item \textbf{Step 1}: This step involves the elimination of unknown
	coefficients from the governing systems, resulting in a boundary value
	problem for a system of nonlinear PDEs. In the case of time dependent data,
	Volterra-like integrals with respect to $t$ are also involved in this system.
	
	\item \textbf{Step 2}: This step focuses on the numerical solution of the
	problem formulated in Step 1. The solution of this problem provides the
	solution for the original CIP. \ To solve that problem, a least squares
	weighted Tikhonov-like functional is constructed. The key element of this
	functional is the presence in it of the Carleman Weight Function (CWF) in
	it, which is not the case of the classical Tikhonov regularization
	functional \cite{T}. The CWF is the function, which is used as the weight in
	the Carleman estimate for the corresponding PDE operator. Furthermore, it is
	demonstrated numerically in many of our publications, such as, e.g. \cite%
	{KTR1,KTR2} as well as in Test 1 in section 8 of this paper that the absence
	of the CWF significantly deteriorates the solution. The key result of the
	convergence analysis of any version of the convexification method is a
	theorem, which claims that, for a proper choice of parameters, that
	functional is strongly convex on a convex set of an arbitrary but fixed
	diameter $d>0$. Next, it is derived from this result that the gradient
	descent method of the minimization of that functional converges to the true
	solution of the CIP\ starting from an arbitrary point of that set, as long
	as the level of noise in the data tends to zero. Since smallness conditions
	are not imposed on $d$, then we call this property \textquotedblleft global
	convergence".
\end{enumerate}

As a by-product, we obtain uniqueness theorem for our CIP. Uniqueness for
similar CIPs for one parabolic PDE was previously proven in \cite%
{Klib92,Ksurvey}, \cite[Theorem 1.10.7]{BK}, \cite[Theorem 3.4.3]{KL}. We
also refer to the Lipschitz stability estimate of \cite{Yam}, which implies
uniqueness. All these works use the framework of the paper \cite{BukhKlib},
in which the tool of Carleman estimates was introduced in the theory of CIPs
for the first time. The convexification method is a numerical development of
the idea of \cite{BukhKlib}, which was originally formulated for a purely
theoretical purpose.

\textbf{Remark 1.1}. \emph{Since the CIP, which we consider, is a quite
	complicated one, then we simplify the presentation via assuming that our
	domain of interest }$\Omega $\emph{\ is a rectangle and all functions we
	work with are sufficiently smooth. It is well known that smoothness
	assumptions are not of a significant concern in the theory of CIPs, see,
	e.g. \cite{Nov1,Nov2}, \cite[Theorem 4.1]{Rom}.}

All functions considered below are real valued ones. In section 2 we pose
both forward and inverse problems. In section 3 we carry out the above Step
1. In section 4 we formulate the Carleman estimate we work with. In \
section 5 we carry out the above Step 2: we construct that functional and
formulate theorems of its global convergence analysis. We prove those
theorems in sections 6 and 7. In section 8 we present our numerical
experiments. Summary of results is given in section 9.

\section{Statements of Forward and Inverse Problems}

\label{sec:2}

Let $N>1$ be an integer and $B$ be a Banach space with its norm $\left\Vert
\cdot \right\Vert _{B}$. Denote%
\begin{equation*}
	\left.
	\begin{array}{c}
		B_{N}=B\times B\times ...\times B,\text{ }N\text{ times,} \\
		\left\Vert u\right\Vert _{B_{N}}=\left( \sum\limits_{k=1}^{N}\left\Vert
		u_{k}\right\Vert _{B}^{2}\right) ^{1/2},\text{ }\forall u=\left(
		u_{1},...,u_{N}\right) \in B_{N}.%
	\end{array}
	\right.
\end{equation*}

Let $a,b,A>0$ be some numbers. Let the number $\eta \in \left( 0,1\right) .$
We set the domain $\Omega \subset \mathbb{R}^{2}$ as a rectangle,
\begin{equation}
	\left.
	\begin{array}{c}
		\Omega =\left\{ \mathbf{x}=\left( x,y\right) :a<x<b,\left\vert y\right\vert
		<A\right\} ,\mbox{ } \\
		\Gamma =\partial \Omega \cap \left\{ x=b\right\} , \\
		Q_{T}=\Omega \times \left( 0,T\right) ,\mbox{ }S_{T}=\partial \Omega \times
		\left( 0,T\right) ,\mbox{ }\Gamma _{T}=\Gamma \times \left( 0,T\right) , \\
		Q_{\eta T}=\Omega \times \left( \left( 1-\eta \right) T/2,\left( 1+\eta
		\right) T/2\right) \subset Q_{T}.%
	\end{array}
	\right.  \label{2.1}
\end{equation}

We now follow notations of \cite{Lee}. Shifting time as $t^{\prime }=t-t_{1}$
and still keeping the same notation for time, we assume that $t\in \left(
0,T\right) .$ Let $\rho _{S}\left( \mathbf{x},t\right) ,\rho _{I}\left(
\mathbf{x},t\right) $ and $\rho _{R}\left( \mathbf{x},t\right) $ are S,I and
R populations respectively. Recall that the 2D vector functions $q_{S}\left(
\mathbf{x}\right) ,q_{I}\left( \mathbf{x}\right) $,$q_{R}\left( \mathbf{x}%
\right) $ are velocities of movements of these populations. Let $\beta
\left( \mathbf{x}\right) $ and $\gamma \left( \mathbf{x}\right) $ be the
infection and recovery rates respectively. The initial boundary value
problem for the spatial SIR model has the form \cite[formulas (2.1)]{Lee}:
\begin{equation}
	\partial _{t}\rho _{S}-\frac{\eta _{S}^{2}}{2}\Delta \rho _{S}+\text{div}
	\left( \rho _{S}q_{S}\right) +\beta \left( \mathbf{x}\right) \rho _{S}\rho
	_{I}=0,\text{ }\left( \mathbf{x},t\right) \in Q_{T},  \label{2.2}
\end{equation}%
\begin{equation}
	\partial _{t}\rho _{I}-\frac{\eta _{I}^{2}}{2}\Delta \rho _{I}+\text{div}
	\left( \rho _{I}q_{I}\right) -\beta \left( \mathbf{x}\right) \rho _{S}\rho
	_{I}=0,\text{ }\left( \mathbf{x},t\right) \in Q_{T},  \label{2.3}
\end{equation}%
\begin{equation}
	\partial _{t}\rho _{R}-\frac{\eta _{R}^{2}}{2}\Delta \rho _{R}+\text{div}
	\left( \rho _{R}q_{R}\right) -\gamma \left( \mathbf{x}\right) \rho _{I}=0,
	\text{ }\left( \mathbf{x},t\right) \in Q_{T},  \label{2.4}
\end{equation}%
\begin{equation}
	\partial _{n}\rho _{S}\mathbf{\mid }_{S_{T}}=g_{1}\left( \mathbf{x},t\right)
	,\text{ }\partial _{n}\rho _{I}\mathbf{\mid }_{S_{T}}=g_{2}\left( \mathbf{x}
	,t\right) ,\text{ }\partial _{n}\rho _{R}\mathbf{\mid }_{S_{T}}=g_{3}\left(
	\mathbf{x},t\right) .  \label{2.5}
\end{equation}%
\begin{equation}
	\rho _{S}\left( \mathbf{x},0\right) =\rho _{S}^{0}\left( \mathbf{x}\right) ,
	\text{ }\rho _{I}\left( \mathbf{x},0\right) =\rho _{I}^{0}\left( \mathbf{x}
	\right) ,\text{ }\rho _{R}\left( \mathbf{x},0\right) =\rho _{R}^{0}\left(
	\mathbf{x}\right) ,  \label{2.6}
\end{equation}%
where $\partial _{n}$ is the normal derivative. This is our forward problem.
Functions $g_{1},g_{2},g_{3}$ are fluxes of those populations through the
boundary, and they are set to zero in \cite{Lee}. However, our theory can
consider a more general case of non-zero fluxes. Here, $\eta _{S}^{2},\eta
_{R}^{2},\eta _{R}^{2}>0$ are constant viscosity terms, which are added in
\cite{Lee} to regularize the system. To simplify the presentation, we assume
below without any loss of generality that%
\begin{equation}
	\frac{\eta _{S}^{2}}{2}=\frac{\eta _{I}^{2}}{2}=\frac{\eta _{R}^{2}}{2}=c>0,
	\label{2.7}
\end{equation}%
where $c$ is a number, although our theory can be easily extended to the
case when condition (\ref{2.7}) is not imposed. Suppose for a moment that
the boundary $\partial \Omega $ is sufficiently smooth, so as functions $%
\beta \left( \mathbf{x}\right) ,\gamma \left( \mathbf{x}\right) ,q_{S}\left(
\mathbf{x}\right) ,q_{I}\left( \mathbf{x}\right) $,$q_{R}\left( \mathbf{x}%
\right) $. Also, assume for a moment that Neumann boundary conditions (\ref%
{2.5}) are replaced with the Dirichlet boundary conditions. Then Theorem 7.1
of Chapter 7 of \cite{Lad} as well as some other techniques developed in
this book apply that, under some additional conditions, including an upper
bound on $T$ (because of the nonlinearity of this system), the so changed
initial boundary value problem has unique solution $\left( \rho _{S},\rho
_{I},\rho _{R}\right) \left( \mathbf{x}\text{,}t\right) ,$ which is
sufficiently smooth in $\overline{Q}_{T}.$ We, however, assume below that
the forward problem (\ref{2.2})-(\ref{2.6}) has unique solution
\begin{equation}
	\rho _{S}\left( \mathbf{x},t\right) ,\rho _{I}\left( \mathbf{x},t\right)
	,\rho _{R}\left( \mathbf{x},t\right) \in C^{6,3}\left( \overline{Q}
	_{T}\right) ,  \label{2.8}
\end{equation}%
see Remark 1.1.

We now pose our Coefficient Inverse Problem. Without a loss of generality,
it is convenient to assume that measurements of spatial distributions of SIR
populations inside of the affected city are conducted at the moment of time $%
t_{0}=T/2.$

\textbf{Coefficient Inverse Problem (CIP).} Let
\begin{equation}
	\rho _{S}\left( \mathbf{x},\frac{T}{2}\right) =p_{1}\left( \mathbf{x}\right)
	,\text{ }\rho _{I}\left( \mathbf{x},\frac{T}{2}\right) =p_{2}\left( \mathbf{%
		x }\right) ,\text{ }\rho _{R}\left( \mathbf{x},\frac{T}{2}\right)
	=p_{3}\left( \mathbf{x}\right) ,  \label{2.9}
\end{equation}%
\begin{equation}
	\rho _{S}\mid _{\Gamma _{T}}=f_{1}\left( y,t\right) ,\rho _{I}\mid _{\Gamma
		_{T}}=f_{2}\left( y,t\right) ,\rho _{S}\mid _{\Gamma _{T}}=f_{1}\left(
	y,t\right) ,\left( y,t\right) \in \Gamma _{T}.  \label{2.10}
\end{equation}%
\emph{Assume that functions in the right hand sides of (\ref{2.9}) and (\ref%
	{2.10}) are known. However, let initial conditions (\ref{2.6}) be unknown.
	Assuming the smoothness conditions (\ref{1.0}) and (\ref{2.8}), find the
	vector function} $\Phi \left( \mathbf{x},t\right) $ \emph{in (\ref{1.1}),
	which is rewritten now as:}%
\begin{equation}
	\Phi \left( \mathbf{x},t\right) =\left( \beta \left( \mathbf{x}\right)
	,\gamma \left( \mathbf{x}\right) ,\rho _{S}\left( \mathbf{x},t\right) ,\rho
	_{I}\left( \mathbf{x},t\right) ,\rho _{R}\left( \mathbf{x},t\right) \right)
	, \text{ }\left( \mathbf{x},t\right) \in Q_{T}.  \label{2.11}
\end{equation}

Thus, conditions (\ref{2.5}) and (\ref{2.10}) mean that we measure fluxes of
SIR populations at the entire boundary of the city, whereas the numbers of
these populations are measured only at the part $\Gamma \subset \partial
\Omega $ of the boundary, see (\ref{2.1}). Hence, (\ref{2.5}) and (\ref{2.10}%
) are incomplete lateral Cauchy data. Conditions (\ref{2.9}) mean that we
measure spatial distributions of SIR populations inside that city at a
single moment of time $t=T/2.$

\section{Transformation}

\label{sec:3}

Following the first step mentioned in section 1, we now transform our CIP in
a boundary value problem for a system of four coupled nonlinear integral
differential equations containing Volterra integrals with respect to $t$.
But this system should not contain unknown coefficients $\beta \left(
\mathbf{x}\right) $ and $\gamma \left( \mathbf{x}\right) .$Assume that in (%
\ref{2.9})
\begin{equation}
	\left\vert p_{1}\left( \mathbf{x}\right) \right\vert ,\left\vert p_{2}\left(
	\mathbf{x}\right) \right\vert \geq \kappa ,\text{ }\mathbf{x}\in \overline{%
		\Omega },  \label{3.1}
\end{equation}%
where $\kappa >0$ is a number. Setting in (\ref{2.2}) and (\ref{2.4}) $t=T/2$
and using (\ref{2.7}), (\ref{2.9}) and (\ref{3.1}), we obtain%
\begin{equation}
	\beta \left( \mathbf{x}\right) =-\frac{1}{\left( p_{1}p_{2}\right) \left(
		\mathbf{x}\right) }\partial _{t}\rho _{S}\left( \mathbf{x},\frac{T}{2}
	\right) +\frac{1}{\left( p_{1}p_{2}\right) \left( \mathbf{x}\right) }\left[
	c\Delta p_{1}\left( \mathbf{x}\right) -\text{div}\left( p_{1}q_{S}\right)
	\left( \mathbf{x}\right) \right] ,  \label{3.2}
\end{equation}%
\begin{equation}
	\gamma \left( \mathbf{x}\right) =\frac{1}{p_{2}\left( \mathbf{x}\right) }
	\partial _{t}\rho _{R}\left( \mathbf{x},\frac{T}{2}\right) -\frac{1}{
		p_{2}\left( \mathbf{x}\right) }\left[ c\Delta p_{3}\left( \mathbf{x}\right)
	- \text{div}\left( p_{3}q_{R}\right) \left( \mathbf{x}\right) \right] .
	\label{3.3}
\end{equation}

Introduce new functions $v_{i}\left( \mathbf{x,}t\right) ,$ $i=1,2,3,$
\begin{equation}
	\left.
	\begin{array}{c}
		v_{1}\left( \mathbf{x,}t\right) =\partial _{t}\rho _{S}\left( \mathbf{x,}
		t\right) ,v_{2}\left( \mathbf{x,}t\right) =\partial _{t}\rho _{I}\left(
		\mathbf{x,}t\right) ,v_{3}\left( \mathbf{x,}t\right) =\partial _{t}\rho
		_{R}\left( \mathbf{x,}t\right) , \\
		V_{1}\left( \mathbf{x,}t\right) =\left( v_{1},v_{2},v_{3}\right) ^{T}\left(
		\mathbf{x,}t\right) .%
	\end{array}
	\right.  \label{3.4}
\end{equation}%
Then by (\ref{2.9})%
\begin{equation}
	\left.
	\begin{array}{c}
		\rho _{S}\left( \mathbf{x,}t\right) =\int\limits_{T/2}^{t}v_{1}\left(
		\mathbf{x,}\tau \right) d\tau +p_{1}\left( \mathbf{x}\right) ,\text{ }\left(
		\mathbf{x,}t\right) \in Q_{T}, \\
		\rho _{I}\left( \mathbf{x,}t\right) =\int\limits_{T/2}^{t}v_{2}\left(
		\mathbf{x,}\tau \right) d\tau +p_{2}\left( \mathbf{x}\right) ,\text{ }\left(
		\mathbf{x,}t\right) \in Q_{T}, \\
		\rho _{R}\left( \mathbf{x,}t\right) =\int\limits_{T/2}^{t}v_{2}\left(
		\mathbf{x,}\tau \right) d\tau +p_{3}\left( \mathbf{x}\right) ,\text{ }\left(
		\mathbf{x,}t\right) \in Q_{T}.%
	\end{array}
	\right.  \label{3.5}
\end{equation}%
Since by (\ref{3.4})%
\begin{equation*}
	\left.
	\begin{array}{c}
		\partial _{t}\rho _{S}\left( \mathbf{x,}T/2\right) =v_{1}\left( \mathbf{x,}
		T/2\right) =v_{1}\left( \mathbf{x,}t\right) -\int\limits_{T/2}^{t}\partial
		_{t}v_{1}\left( \mathbf{x,}\tau \right) d\tau , \\
		\partial _{t}\rho _{R}\left( \mathbf{x,}T/2\right) =v_{3}\left( \mathbf{x,}
		T/2\right) =v_{3}\left( \mathbf{x,}t\right) -\int\limits_{T/2}^{t}\partial
		_{t}v_{3}\left( \mathbf{x,}\tau \right) d\tau ,%
	\end{array}
	\right.
\end{equation*}%
then formulas (\ref{3.2}) and (\ref{3.3}) become:%
\begin{equation}
	\left.
	\begin{array}{c}
		\beta \left( \mathbf{x}\right) =\left( v_{1}\left( \mathbf{x,}t\right)
		-\int\limits_{T/2}^{t}\partial _{t}v_{1}\left( \mathbf{x,}\tau \right) d\tau
		\right) r_{1}\left( \mathbf{x}\right) +r_{2}\left( \mathbf{x}\right) , \\
		\gamma \left( \mathbf{x}\right) =\left( v_{3}\left( \mathbf{x,}t\right)
		-\int\limits_{T/2}^{t}\partial _{t}v_{3}\left( \mathbf{x,}\tau \right) d\tau
		\right) r_{3}\left( \mathbf{x}\right) +r_{4}\left( \mathbf{x}\right) , \\
		r_{1}\left( \mathbf{x}\right) =-\left[ \left( p_{1}p_{2}\right) \left(
		\mathbf{x}\right) \right] ^{-1}, \\
		r_{2}\left( \mathbf{x}\right) =-r_{1}\left( \mathbf{x}\right) \left[ c\Delta
		p_{1}\left( \mathbf{x}\right) -\text{div}\left( p_{1}q_{S}\right) \left(
		\mathbf{x}\right) \right] , \\
		r_{3}\left( \mathbf{x}\right) =1/p_{2}\left( \mathbf{x}\right) ,\text{ } \\
		r_{4}\left( \mathbf{x}\right) =-r_{3}\left( \mathbf{x}\right) \left[ c\Delta
		p_{3}\left( \mathbf{x}\right) -\text{div}\left( p_{3}q_{R}\right) \left(
		\mathbf{x}\right) \right] .%
	\end{array}
	\right.  \label{3.6}
\end{equation}%
Differentiate (\ref{2.2})-(\ref{2.5}) and (\ref{2.10}) with respect to $t$
and use (\ref{3.4})-(\ref{3.6}). We obtain the following system of three
nonlinear integral differential equations with incomplete lateral Cauchy
data:%
\begin{equation}
	\left.
	\begin{array}{c}
		\partial _{t}V_{1}-c\Delta V_{1}-P_{1}\left( V_{1},\nabla
		V_{1},\int\limits_{T/2}^{t}V_{1}\left( \mathbf{x},\tau \right) d\tau
		,\int\limits_{T/2}^{t}\partial _{t}V_{1}\left( \mathbf{x},\tau \right) d\tau
		,\mathbf{x}\right) =0\text{ in }Q_{T}, \\
		\partial _{n}V_{1}\mid _{S_{T}}=\left( \partial _{t}g_{1},\partial
		_{t}g_{2},\partial _{t}g_{3}\right) ^{T}\left( \mathbf{x},t\right) ,\text{ }
		V_{1}\mid _{\Gamma _{T}}=\left( \partial _{t}f_{1},\partial
		_{t}f_{2},\partial _{t}f_{3}\right) ^{T}\left( y,t\right) ,%
	\end{array}
	\right.  \label{3.7}
\end{equation}%
where the 3D vector function $P_{1}$ is infinitely many times differentiable
with respect to all of its variables, except of $\mathbf{x}$ and it is
continuous with respect to $\mathbf{x}$.

In terms of Carleman estimates, an inconvenient feature of system (\ref{3.7}%
) is the presence of the integral containing the time derivative $\partial
_{t}V_{1}$. Indeed, since the term $\partial _{t}V_{1}$ is involved in the
principal part $\partial _{t}-c\Delta $ of the operator in (\ref{3.7}), then
it is unclear how to absorb this term using the Carleman estimate of section
4.

Therefore, we differentiate (\ref{3.7}) with respect to $t$ and denote
\begin{equation}
	V_{2}\left( \mathbf{x},t\right) =\partial _{t}V_{1}\left( \mathbf{x}
	,t\right) =\left( w_{1},w_{2},w_{3}\right) ^{T}\left( \mathbf{x},t\right) .
	\label{3.8}
\end{equation}%
We obtain%
\begin{equation}
	\left.
	\begin{array}{c}
		\partial _{t}V_{2}-c\Delta V_{2}- \\
		-P_{2}\left( V_{1},\nabla V_{1},V_{2},\nabla
		V_{2}\int\limits_{T/2}^{t}V_{1}\left( \mathbf{x},\tau \right) d\tau
		,\int\limits_{T/2}^{t}V_{2}\left( \mathbf{x},\tau \right) d\tau ,\mathbf{x}
		\right) =0\text{ in }Q_{T}, \\
		\partial _{n}V_{2}\mid _{S_{T}}=\left( \partial _{t}^{2}g_{1},\partial
		_{t}^{2}g_{2},\partial _{t}^{2}g_{3}\right) ^{T}\left( \mathbf{x},t\right) ,
		\text{ }V_{2}\mid _{\Gamma _{T}}=\left( \partial _{t}^{2}f_{1},\partial
		_{t}^{2}f_{2},\partial _{t}^{2}f_{3}\right) ^{T}\left( \mathbf{x},t\right) ,%
	\end{array}
	\right.  \label{3.9}
\end{equation}%
where the 3D vector function $P_{2}$ has the same smoothness properties as
the ones of $P_{1}$. Using (\ref{3.4}) and (\ref{3.7})-(\ref{3.9}), consider
now the 6D vector functions
\begin{equation}
	\left.
	\begin{array}{c}
		W\left( \mathbf{x},t\right) =\left(
		v_{1},v_{2},v_{3},w_{1},w_{2},w_{3}\right) ^{T}\left( \mathbf{x,}t\right)
		\in C_{6}^{2,1}\left( \overline{Q}_{T}\right) ,\text{ } \\
		G\left( \mathbf{x,}t\right) =\left( \partial _{t}g_{1},\partial
		_{t}g_{2},\partial _{t}g_{3},\partial _{t}^{2}g_{1},\partial
		_{t}^{2}g_{2},\partial _{t}^{2}g_{3}\right) ,\text{ }\left( \mathbf{x,}
		t\right) \in S_{T}, \\
		F\left( \mathbf{x,}t\right) =\left( \partial _{t}f_{1},\partial
		_{t}f_{2},\partial _{t}f_{3},\partial _{t}^{2}f_{1},\partial
		_{t}^{2}f_{2},\partial _{t}^{2}f_{3}\right) ,\text{ }\left( \mathbf{x,}
		t\right) \in \Gamma _{T}.%
	\end{array}
	\right.  \label{3.10}
\end{equation}

Using (\ref{3.7})-(\ref{3.10}), we obtain the following boundary value
problem with incomplete lateral Cauchy data for the $6\times 6$ system of
coupled nonlinear integral differential equations:
\begin{equation}
	\left.
	\begin{array}{c}
		L\left( W\right) =W_{t}-c\Delta W-P\left( W,\nabla
		W,\int\limits_{T/2}^{t}W\left( \mathbf{x},\tau \right) d\tau ,\mathbf{x}
		\right) =0\text{ in }Q_{T}, \\
		\partial _{n}W\mid _{S_{T}}=G\left( \mathbf{x},t\right) ,\text{ }W\mid
		_{\Gamma _{T}}=F\left( \mathbf{x},t\right) ,%
	\end{array}
	\right.  \label{3.11}
\end{equation}%
where the 6D vector function $P$ is infinitely many times differentiable
with respect to all of its variables, except of $\mathbf{x,}$ and it is
continuous with respect to $\mathbf{x}\in \overline{\Omega }$.

Suppose that the solution $W\left( \mathbf{x},t\right) \in C_{6}^{2,1}\left(
\overline{Q}_{T}\right) $ of problem (\ref{3.11}) is computed. Then the
first line of (\ref{3.10}) implies that the unknown infection and recovery
rates $\beta \left( \mathbf{x}\right) $ and $\gamma \left( \mathbf{x}\right)
$ are found immediately via (\ref{3.6}). Finally, our main target, the
spatiotemporal distributions of S,I and R populations, are found immediately
via (\ref{3.5}). Therefore, we focus below on the numerical solution of
problem (\ref{3.11}).

One might also raise the question about stable computations of first and
second $t-$derivatives of boundary data, since these data are usually noisy.
This research group is experienced in the regularization of the numerical
computations of both first \cite{KLZpar} and second \cite{MFG8} derivatives
of noisy data via cubic spline interpolation, also, see section 8.

\section{Carleman Estimate}

\label{sec:4}

There are several CWFs, which can be used for Carleman estimates for
parabolic operators, see, e.g. \cite[Theorem 2.3.1]{KL}, \cite[Chapter 4, \S %
1]{LRS}, \cite{Yam}. Those CWFs depend on two large parameters. However,
since they changes too rapidly, then this is inconvenient for the
convexification, in which a CWF is involved in the numerical scheme. Hence,
we use a CWF, which depends only on a single large parameter. \ More
precisely, our CWF is%
\begin{equation}
	\varphi _{\lambda }\left( \mathbf{x},t\right) =\exp \left[ 2\lambda \left(
	x^{2}-\left( t-T/2\right) ^{2}\right) \right] ,  \label{4.1}
\end{equation}%
where $\lambda \geq 1$ is a large parameter. It follows from (\ref{2.1}) and
(\ref{4.1}) that
\begin{equation}
	\varphi _{\lambda }\left( \mathbf{x},t\right) \leq \exp \left[ 2\lambda
	\left( b^{2}-\left( t-T/2\right) ^{2}\right) \right] \text{ in }Q_{T},
	\label{4.2}
\end{equation}%
\begin{equation}
	\varphi _{\lambda }\left( \mathbf{x},t\right) \geq \exp \left[ -2\lambda
	\left( t-T/2\right) ^{2}\right] \text{ in }Q_{T}.  \label{4.20}
\end{equation}%
Introduce the subspace $H_{0}^{2,1}\left( Q_{T}\right) $ of the space $%
H^{2,1}\left( Q_{T}\right) $ as:%
\begin{equation}
	H_{0}^{2,1}\left( Q_{T}\right) =\left\{ u\in H^{2,1}\left( Q_{T}\right)
	:\partial _{n}u\mid _{S_{T}}=0,\text{ }u\mid _{\Gamma _{T}}=0\right\} .
	\label{4.3}
\end{equation}

\textbf{Theorem 4.1} (Carleman estimate for the operator $\partial
_{t}-c\Delta $ \cite{KLZpar}). \emph{Temporary denote }$\mathbf{x}=\left(
x,y\right) =\left( x_{1},x_{2}\right) .$\emph{\ Let }$c>0$\emph{\ be the
	number in (\ref{2.7}). There exists a sufficiently large number }$\lambda
_{0}=\lambda _{0}\left( c,Q_{T}\right) \geq 1$\emph{\ and a number }$%
C=C\left( c,Q_{T}\right) >0,$\emph{\ both depending only on listed
	parameters, such that the following Carleman estimate holds:}%
\begin{equation*}
	\int\limits_{Q_{T}}\left( u_{t}-c\Delta u\right) ^{2}\varphi _{\lambda }d
	\mathbf{x}dt\geq \frac{C}{\lambda }\int\limits_{Q_{T}}u_{t}^{2}\varphi
	_{\lambda }d\mathbf{x}dt+\frac{C}{\lambda }\sum\limits_{i,j=1}^{2}\int
	\limits_{Q_{T}}u_{x_{i}x_{j}}^{2}\varphi _{\lambda }d\mathbf{x}dt+
\end{equation*}%
\begin{equation*}
	+C\int\limits_{Q_{T}}\left( \lambda \left( \nabla u\right) ^{2}+\lambda
	^{3}u^{2}\right) \varphi _{\lambda }d\mathbf{x}dt-
\end{equation*}%
\begin{equation*}
	-C\left( \left\Vert u\left( \mathbf{x},T\right) \right\Vert _{H^{1}\left(
		\Omega \right) }^{2}+\left\Vert u\left( \mathbf{x},0\right) \right\Vert
	_{H^{1}\left( \Omega \right) }^{2}\right) \lambda ^{2}\exp \left( -2\lambda
	\left( \frac{T^{2}}{4}-b^{2}\right) \right) ,
\end{equation*}%
\begin{equation*}
	\forall u\in H_{0}^{2,1}\left( Q_{T}\right) ,\text{ }\forall \lambda \geq
	\lambda _{0}.
\end{equation*}

This Carleman estimate is proven in \cite{KLZpar} for the case when boundary
conditions (\ref{4.3}) are replaced with $u\mid _{S_{T}}=0,$ $u_{x}\mid
_{\Gamma _{T}}=0.$ However, it easily follows from that proof that the same
result is valid for boundary conditions (\ref{4.3}): we refer to formulas
(9.4) and (9.7) of \cite{KLZpar}.

We also need an estimate of the Volterra integral, in which the CWF (\ref%
{4.1}) is involved.

\textbf{Lemma 4.1 }(\cite{Ksurvey}, \cite[Lemma 3.1.1]{KL}). \emph{The
	following estimate is valid:}%
\begin{equation*}
	\int\limits_{Q_{T}}\left( \int\limits_{T/2}^{t}f\left( \mathbf{x},\tau
	\right) d\tau \right) ^{2}\varphi _{\lambda }\left( \mathbf{x},t\right) d
	\mathbf{x}dt\leq \frac{C}{\lambda }\int\limits_{Q_{T}}f^{2}\left( \mathbf{x}
	,t\right) \varphi _{\lambda }\left( \mathbf{x},t\right) d\mathbf{x}dt,\text{
	}\forall \lambda >0,\forall f\in L_{2}\left( Q_{T}\right) ,
\end{equation*}%
\emph{where the number }$C=C\left( Q_{T}\right) >0$\emph{\ depends only on
	the domain }$Q_{T}.$

\section{The Convexification Functional}

\label{sec:5}

\subsection{The functional}

\label{sec:5.1}

For our construction we need that the solution of problem (\ref{3.11}) $W\in
C_{6}^{2,1}\left( \overline{Q}_{T}\right) .$ To apply the Sobolev embedding
theorem, we require a little bit more: that $W\in C_{6}^{2}\left( \overline{Q%
}_{T}\right) .$ Since $Q_{T}\subset \mathbb{R}^{3},$ then Sobolev embedding
theorem implies that
\begin{equation}
	H^{4}\left( Q_{T}\right) \subset C^{2}\left( \overline{Q}_{T}\right) ,\text{
	}\left\Vert u\right\Vert _{C^{2}\left( \overline{Q}_{T}\right) }\leq
	C\left\Vert u\right\Vert _{H^{4}\left( Q_{T}\right) },\text{ }\forall u\in
	H^{4}\left( Q_{T}\right)  \label{5.1}
\end{equation}%
with a certain constant $C=C\left( Q_{T}\right) >0$ depending only on the
domain $Q_{T}.$

Let $R>0$ be an arbitrary number. Keeping in mind the lateral Cauchy data in
(\ref{3.11}), define the sets $B\left( R\right) ,B_{0}\left( R\right) $ of
6D vector functions as:%
\begin{equation}
	B\left( R\right) =\left\{ U\in H_{6}^{4}\left( Q_{T}\right) :\partial
	_{n}U\mid _{S_{T}}=G\left( \mathbf{x},t\right) ,\text{ }U\mid _{\Gamma
		_{T}}=F\left( y,t\right) ,\text{ }\left\Vert U\right\Vert _{H_{6}^{4}\left(
		Q_{T}\right) }<R\right\} .  \label{5.2}
\end{equation}%
By (\ref{5.1}) and (\ref{5.2})%
\begin{equation}
	B\left( R\right) \subset C_{6}^{2}\left( \overline{Q}_{T}\right) ,\text{ }
	\left\Vert U\right\Vert _{C_{6}^{2}\left( \overline{Q}_{T}\right) }\leq CR,
	\text{ }\forall U\in \overline{B\left( R\right) }.  \label{5.3}
\end{equation}

Let $\xi \in \left( 0,1\right) $ be the regularization parameter and $%
L\left( W\right) $ be the nonlinear integral differential operator in (\ref%
{3.11}). We define our regularization convexification functional $J_{\lambda
	,\xi }:\overline{B\left( R\right) }\rightarrow \mathbb{R}$ as:%
\begin{equation}
	J_{\lambda ,\xi }\left( W\right) =e^{-2\lambda
		b^{2}}\int\limits_{Q_{T}}\left( L\left( W\right) \right) ^{2}\varphi
	_{\lambda }d\mathbf{x}dt+\xi \left\Vert W\right\Vert _{H_{6}^{4}\left(
		Q_{T}\right) }^{2}.  \label{5.4}
\end{equation}%
The multiplier $e^{-2\lambda b^{2}}$ in (\ref{5.4}) is the balancing factor
between two terms in the right hand side of (\ref{5.4}), see (\ref{4.2}). To
find the desired approximate solution of problem (\ref{3.11}), we consider
below the following problem:

\textbf{Minimization Problem}. \emph{Minimize the functional }$J_{\lambda
	,\xi }\left( W\right) $\emph{\ on the set }$\overline{B\left( R\right) }.$

\subsection{Strong convexity of the functional $J_{\protect\lambda ,\protect%
		\xi }\left( W\right) $ on the set $\overline{B\left( R\right) }$}

\label{sec:5.2}

Define the subspace $H_{6,0}^{4}\left( Q_{T}\right) $ of the space $%
H_{6}^{4}\left( Q_{T}\right) $ as:%
\begin{equation}
	H_{6,0}^{4}\left( Q_{T}\right) =\left\{ U\in H_{6}^{4}\left( Q_{T}\right)
	:\partial _{n}U\mid _{S_{T}}=0,\text{ }U\mid _{\Gamma _{T}}=0\right\} .
	\label{5.5}
\end{equation}

\textbf{Theorem 5.1} (strong convexity of $J_{\lambda ,\beta }\left(
W\right) $). \emph{Assume that inequalities (\ref{3.1}) hold. The following
	statements are true:}

\emph{1. The functional }$J_{\lambda ,\beta }\left( W\right) $\emph{\ has
	the Fr\'{e}chet derivative }$J_{\lambda ,\xi }^{\prime }\left( W\right) \in
H_{6,0}^{4}\left( Q_{T}\right) $\emph{\ at every point }$W\in \overline{
	B\left( R\right) }$\emph{\ and for every pair of numbers }$\lambda >0,\xi
>0. $\emph{\ This derivative is Lipschitz continuous in }$\overline{B\left(
	R\right) }.$\emph{\ The latter means that there exists a number }$D=D\left(
\kappa ,\lambda ,c,\xi ,R,Q_{T}\right) $\emph{\ such that the following
	inequality holds:}%
\begin{equation}
	\left\Vert J_{\lambda ,\xi }^{\prime }\left( W_{1}\right) -J_{\lambda ,\xi
	}^{\prime }\left( W_{2}\right) \right\Vert _{H_{6}^{4}\left( Q_{T}\right)
	}\leq D\left\Vert W_{1}-W_{2}\right\Vert _{H_{6}^{4}\left( Q_{T}\right) },
	\text{ }\forall W_{1},W_{2}\in \overline{B\left( R\right) }.  \label{5.6}
\end{equation}

\emph{2. Temporary denote again }$\mathbf{x}=\left( x,y\right) =\left(
x_{1},x_{2}\right) .$\emph{\ Let }$\lambda _{0}\geq \lambda _{0}\left(
c,Q_{T}\right) \geq 1$\emph{\ be the number of Theorem 4.1. There exists a
	sufficiently large number }$\lambda _{1}=\lambda _{1}\left( \kappa
,c,R,Q_{T}\right) \geq \lambda _{0}$\emph{\ and a number }$C_{1}=C_{1}\left(
\kappa ,c,R,Q_{T}\right) >0$\emph{, \ both depending only on listed
	parameters, such that if }$\lambda \geq \lambda _{1}$\emph{\ and if the
	regularization parameter }$\xi $\emph{\ is such that}
\begin{equation}
	\frac{\xi }{2}\in \left[ \exp \left( -\lambda \frac{T^{2}}{4}\right) ,\frac{%
		1 }{2}\right) ,  \label{5.7}
\end{equation}%
\emph{then the functional }$J_{\lambda ,\xi }\left( W\right) $\emph{\ is
	strongly convex on }$\overline{B\left( R\right) }.$\emph{\ More precisely,
	the following inequality holds:}%
\begin{equation*}
	J_{\lambda ,\xi }\left( W_{2}\right) -J_{\lambda ,\xi }\left( W_{1}\right)
	-J_{\lambda ,\xi }^{\prime }\left( W_{1}\right) \left( W_{2}-W_{1}\right)
	\geq
\end{equation*}%
\begin{equation*}
	\geq \frac{C_{1}}{\lambda }\exp \left( -\lambda \frac{T^{2}}{2}\right)
	\left\Vert \partial _{t}\left( W_{2}-W_{1}\right) \right\Vert _{L_{2}\left(
		Q_{T}\right) }^{2}+
\end{equation*}%
\begin{equation}
	+\frac{C_{1}}{\lambda }\exp \left( -\lambda \frac{T^{2}}{2}\right)
	\sum\limits_{i,j=1}^{2}\left\Vert \left( W_{2}-W_{1}\right)
	_{x_{i}x_{j}}\right\Vert _{L_{2}\left( Q_{T}\right) }^{2}+  \label{5.8}
\end{equation}%
\begin{equation*}
	+C_{1}\exp \left( -\lambda \frac{T^{2}}{2}\right) \left( \lambda \left\Vert
	\nabla \left( W_{2}-W_{1}\right) \right\Vert _{L_{2}\left( Q_{T}\right)
	}^{2}+\lambda ^{3}\left\Vert \nabla \left( W_{2}-W_{1}\right) \right\Vert
	_{L_{2}\left( Q_{T}\right) }^{2}\right) +
\end{equation*}%
\begin{equation*}
	+\frac{\xi }{2}\left\Vert W_{2}-W_{1}\right\Vert _{H^{4}\left( Q_{T}\right)
	}^{2},
\end{equation*}%
\begin{equation*}
	\text{ }\forall W_{1},W_{2}\in \overline{B\left( R\right) },\text{ }\forall
	\lambda \geq \lambda _{1}.
\end{equation*}

\emph{3. For any pair of numbers }$\left( \lambda ,\xi \right) $\emph{\ of
	item 2 there exists unique minimizer }$W_{\min ,\lambda ,\xi }\in \overline{
	B\left( R\right) }$\emph{\ of the functional }$J_{\lambda ,\xi }\left(
W\right) $\emph{\ on the set }$\overline{B\left( R\right) }$\emph{.
	Furthermore, the following inequality holds:}%
\begin{equation}
	J_{\lambda ,\xi }^{\prime }\left( W_{\min ,\lambda ,\xi }\right) \left(
	W-W_{\min ,\lambda ,\xi }\right) \geq 0,\text{ }\forall W\in \overline{
		B\left( R\right) }.  \label{5.9}
\end{equation}

\textbf{Remarks 5.1}:

\begin{enumerate}
	\item \emph{Even though sufficiently large values of }$\lambda $\emph{\ are
		required in this theorem, our past computational experience of, e.g. \cite%
		{KL}-\cite{MFG8} for several versions of the convexification as well as of
		this paper (section 8) demonstrates that }$\lambda \in \left[ 1,5\right] $
	\emph{\ is sufficient. In fact, this is similar with any asymptotic theory.
		Such a theory usually states that if a certain parameter }$X_{1}$\emph{\ is
		sufficiently large, then a certain formula }$X_{2}$\emph{\ is valid with a
		good accuracy. However, in practical computations with specific ranges of
		parameters of specific mathematical models only computational experience can
		tell one which exactly values of }$X_{1}$\emph{\ provide a good accuracy for
		the formula }$X_{2}$\emph{.}
	
	\item \emph{The minimizer }$W_{\min ,\lambda ,\xi }$\emph{\ of the
		functional }$J_{\lambda ,\xi }\left( W\right) $\emph{\ is called in the
		regularization theory \cite{T} the \textquotedblleft regularized solution"
		of the CIP (\ref{2.9}), (\ref{2.10}). It is always important in this theory
		to estimate the accuracy of the regularized solution, and this is done in
		subsection 5.3.}
	
	\item \emph{Below }$C_{1}=C_{1}\left( \kappa ,c,R,Q_{T}\right) >0$\emph{\
		denotes different numbers depending only on listed parameters.}
\end{enumerate}

\subsection{The accuracy of the minimizer}

\label{sec:5.3}

We now estimate the accuracy of the minimizer of a slightly modified
functional (\ref{5.4}) with respect to the noise level $\delta \in \left(
0,1\right) $ in the vector function of the input data $\left( p_{1}\left(
\mathbf{x}\right) ,p_{2}\left( \mathbf{x}\right) ,p_{3}\left( \mathbf{x}%
\right) ,G\left( \mathbf{x},t\right) ,F\left( y,t\right) \right) $ in (\ref%
{2.9}) and (\ref{3.11}). We assume that $\delta $ is sufficiently small. In
accordance with the concept of Tikhonov for ill-posed problems \cite{T}, we
assume the existence of the exact vector function $\Phi ^{\ast }\left(
\mathbf{x},t\right) $ in (\ref{2.11})
\begin{equation*}
	\Phi ^{\ast }\left( \mathbf{x},t\right) =\left( \beta ^{\ast }\left( \mathbf{%
		\ x}\right) ,\gamma ^{\ast }\left( \mathbf{x}\right) ,\rho _{S}^{\ast
	}\left( \mathbf{x},t\right) ,\rho _{I}^{\ast }\left( \mathbf{x},t\right)
	,\rho _{R}^{\ast }\left( \mathbf{x},t\right) \right) ,\text{ }\left( \mathbf{%
		x} ,t\right) \in Q_{T}
\end{equation*}%
with the noiseless data in (\ref{2.5}), (\ref{2.9}) and (\ref{2.10}),
\begin{equation}
	\left.
	\begin{array}{c}
		\rho _{S}^{\ast }\left( \mathbf{x},T/2\right) =p_{1}^{\ast }\left( \mathbf{x}
		\right) ,\text{ }\rho _{I}^{\ast }\left( \mathbf{x},T/2\right) =p_{2}^{\ast
		}\left( \mathbf{x}\right) ,\text{ }\rho _{R}^{\ast }\left( \mathbf{x}
		,T/2\right) =p_{3}^{\ast }\left( \mathbf{x}\right) , \\
		\partial _{n}\rho _{S}^{\ast }\mathbf{\mid }_{S_{T}}=g_{1}^{\ast }\left(
		\mathbf{x},t\right) ,\text{ }\partial _{n}\rho _{I}^{\ast }\mathbf{\mid }
		_{S_{T}}=g_{2}^{\ast }\left( \mathbf{x},t\right) ,\text{ }\partial _{n}\rho
		_{R}^{\ast }\mathbf{\mid }_{S_{T}}=g_{3}^{\ast }\left( \mathbf{x},t\right) ,
		\\
		\rho _{S}^{\ast }\mid _{\Gamma _{T}}=f_{1}^{\ast }\left( y,t\right) ,\text{ }
		\rho _{I}^{\ast }\mid _{\Gamma _{T}}=f_{2}^{\ast }\left( y,t\right) ,\text{ }
		\rho _{S}^{\ast }\mid _{\Gamma _{T}}=f_{1}^{\ast }\left( y,t\right) ,\left(
		y,t\right) \in \Gamma _{T}.%
	\end{array}
	\right.  \label{5.12}
\end{equation}%
Hence, obvious direct analogs of formulas (\ref{3.6}) are valid for
functions $\beta ^{\ast }\left( \mathbf{x}\right) $ and $\gamma ^{\ast
}\left( \mathbf{x}\right) .$ Furthermore, we obtain the direct analog $%
W^{\ast }\left( \mathbf{x},t\right) $ of the vector function $W\left(
\mathbf{x},t\right) $ as well as the direct analog of (\ref{3.11}):
\begin{equation}
	\left.
	\begin{array}{c}
		\widetilde{L}\left( W^{\ast }\right) =W_{t}^{\ast }-c\Delta W^{\ast
		}-P^{\ast }\left( W^{\ast },\nabla W^{\ast },\int\limits_{T/2}^{t}W^{\ast
		}\left( \mathbf{x},\tau \right) d\tau ,\mathbf{x}\right) \\
		=0\text{ in }Q_{T}, \text{ \hspace{5.2 cm}} \\
		\partial _{n}W^{\ast }\mid _{S_{T}}=G^{\ast }\left( \mathbf{x},t\right) ,
		W^{\ast }\mid _{\Gamma _{T}}=F^{\ast }\left( y,t\right) . \text{ \hspace{2 cm%
		}}%
	\end{array}
	\right.  \label{5.13}
\end{equation}%
The difference between operators $P$ and $P^{\ast }$ in (\ref{3.11}) and (%
\ref{5.13}) is that functions $p_{1},p_{2},p_{3}$ in $P$ are replaced in $%
P^{\ast }$ with functions $p_{1}^{\ast },p_{2}^{\ast },p_{3}^{\ast },$ which
are defined in the first line of (\ref{5.12}).

We assume that
\begin{equation}
	W^{\ast }\in B^{\ast }\left( R\right) =\left\{
	\begin{array}{c}
		U\in H_{6}^{4}\left( Q_{T}\right) :\partial _{n}U\mid _{S_{T}}=G^{\ast
		}\left( \mathbf{x},t\right) , \\
		U^{\ast }\mid _{\Gamma _{T}}=F^{\ast }\left( y,t\right) ,\text{ }\left\Vert
		U\right\Vert _{H_{6}^{4}\left( Q_{T}\right) }<R%
	\end{array}
	\right\} .  \label{5.14}
\end{equation}%
By (\ref{2.9}), (\ref{5.1}), (\ref{5.12}) and (\ref{5.14}) functions $%
p_{i},p_{i}^{\ast }\in C^{2}\left( \overline{\Omega }\right) ,i=1,2,3.$
Hence, we assume that
\begin{equation}
	\left\Vert p_{i}-p_{i}^{\ast }\right\Vert _{C^{2}\left( \overline{\Omega }
		\right) }<\delta ,\text{ }i=1,2,3.  \label{5.15}
\end{equation}%
Hence, by (\ref{3.1}) we can assume that
\begin{equation}
	\left\vert p_{1}^{\ast }\left( \mathbf{x}\right) \right\vert ,\left\vert
	p_{2}^{\ast }\left( \mathbf{x}\right) \right\vert \geq \kappa ,\text{ }
	\mathbf{x}\in \Omega .  \label{5.150}
\end{equation}%
We also assume that there exist 6D vector functions $Y\left( \mathbf{x}%
,t\right) $ and $Y^{\ast }\left( \mathbf{x},t\right) $ such that
\begin{equation}
	\left.
	\begin{array}{c}
		Y,Y^{\ast }\in H^{4}\left( Q_{T}\right) ,\text{ }\left\Vert Y\right\Vert
		_{H_{6}^{4}\left( Q_{T}\right) }<R,\text{ }\left\Vert Y^{\ast }\right\Vert
		_{H_{6}^{4}\left( Q_{T}\right) }<R, \\
		\partial _{n}Y\mid _{S_{T}}=G\left( \mathbf{x},t\right) ,\text{ }Y\mid
		_{\Gamma _{T}}=F\left( y,t\right) , \\
		\partial _{n}Y^{\ast }\mid _{S_{T}}=G^{\ast }\left( \mathbf{x},t\right) ,
		\text{ }Y^{\ast }\mid _{\Gamma _{T}}=F^{\ast }\left( y,t\right) .%
	\end{array}
	\right.  \label{5.16}
\end{equation}%
In addition, we assume that%
\begin{equation}
	\left\Vert Y-Y^{\ast }\right\Vert _{C^{2,1}\left( \overline{Q}_{T}\right)
	}<\delta .  \label{5.17}
\end{equation}%
Similarly with (\ref{5.2}), we define the set $B_{0}\left( R\right) $ as%
\begin{equation}
	B_{0}\left( R\right) =\left\{ U\in H_{6}^{4}\left( Q_{T}\right) :\partial
	_{n}U\mid _{S_{T}}=0,\text{ }U\mid _{\Gamma _{T}}=0,\text{ }\left\Vert
	U\right\Vert _{H_{6}^{4}\left( Q_{T}\right) }<R\right\} .  \label{5.18}
\end{equation}%
For every vector function $W\in B\left( R\right) $ as well as for the exact
solution $W^{\ast }\in B^{\ast }\left( R\right) $ we define vector functions
$\widetilde{W}$ and $\widetilde{W^{\ast }}$ as:%
\begin{equation}
	\widetilde{W}=W-Y,\text{ }\widetilde{W^{\ast }}=W^{\ast }-Y^{\ast }.
	\label{5.19}
\end{equation}%
Then triangle inequality, (\ref{5.2}), (\ref{5.14}), the first line of (\ref%
{5.16}), (\ref{5.18}) and (\ref{5.19}) imply
\begin{equation}
	\left.
	\begin{array}{c}
		\widetilde{W},\text{ }\widetilde{W^{\ast }}\in B_{0}\left( 2R\right) , \\
		X+Y\in B\left( 3R\right) ,\text{ }\forall X\in B_{0}\left( 2R\right) .%
	\end{array}
	\right.  \label{5.20}
\end{equation}

Consider another functional $I_{\lambda ,\xi },$
\begin{equation}
	I_{\lambda ,\xi }:B_{0}\left( 2R\right) \rightarrow \mathbb{R},\text{ }
	I_{\lambda ,\xi }\left( X\right) =J_{\lambda ,\xi }\left( X+Y\right) .
	\label{5.21}
\end{equation}%
It follows from (\ref{5.18}), (\ref{5.20}) and (\ref{5.21}) that a direct
analog of Theorem 5.1 is valid for the functional $I_{\lambda ,\xi }\left(
X\right) $. In this case $B\left( R\right) $ should be replaced with $%
B_{0}\left( 2R\right) $ and the number $\lambda _{1}=\lambda _{1}\left(
\kappa ,c,R,Q_{T}\right) $ should be replaced with the number $\lambda _{2},$
where
\begin{equation}
	\lambda _{2}=\lambda _{1}\left( \kappa ,c,3R,Q_{T}\right) \geq \lambda _{1}.
	\label{5.22}
\end{equation}

\textbf{Theorem 5.2} (the accuracy of the minimizer). \emph{Suppose that
	conditions (\ref{3.1}), (\ref{5.14})-(\ref{5.22}) hold. For any }$\lambda
\geq \lambda _{2}$\emph{\ let the regularization parameter }%
\begin{equation}
	\xi =2\exp \left( -\lambda \frac{T^{2}}{4}\right) ,  \label{5.23}
\end{equation}%
\emph{\ (see (\ref{5.7})) and let }$X_{\min ,\lambda ,\xi }\in \overline{
	B_{0}\left( 2R\right) }$\emph{\ be the unique minimizer of the functional }$%
I_{\lambda ,\xi }\left( X\right) $\emph{\ on the set }$\overline{B_{0}\left(
	2R\right) },$\emph{\ the existence of which is guaranteed by Theorem 5.1.
	Let }$\mu _{1},\mu _{2}\in \left( 0,1\right) $\emph{\ be two arbitrary
	numbers such that }$\mu _{2}\in \left( 0,\mu _{1}\right) .$\emph{\ There
	exists a sufficiently small number }$\delta _{0}=\delta _{0}\left( \kappa
,c,R,\mu _{2},Q_{T}\right) \in \left( 0,1\right) $\emph{\ depending only on
	listed parameters such that if }%
\begin{equation}
	\left.
	\begin{array}{c}
		\ln \left( \delta _{0}^{-8/T^{2}}\right) >\lambda _{2},\text{ } \\
		\lambda =\lambda \left( \delta \right) =\ln \left( \delta ^{-8/T^{2}}\right)
		,\text{ }\delta \in \left( 0,\delta _{0}\right) \text{,}%
	\end{array}
	\right.  \label{5.24}
\end{equation}%
\begin{equation}
	\eta ^{2}=\frac{\mu _{1}-\mu _{2}}{2},  \label{5.25}
\end{equation}%
\emph{then}
\begin{equation}
	\left.
	\begin{array}{c}
		\left\Vert \left( X_{\min ,\lambda ,\xi }+Y\right) -W^{\ast }\right\Vert
		_{H^{2,1}\left( Q_{\eta T}\right) }\leq C_{1}\delta ^{1-\mu _{1}}, \\
		\left\Vert \rho _{S,\min ,\lambda ,\xi }-\rho _{S}^{\ast }\right\Vert
		_{H^{2,1}\left( Q_{\eta T}\right) }+\left\Vert \rho _{I,\min ,\lambda ,\xi
		}-\rho _{S}^{\ast }\right\Vert _{H^{2,1}\left( Q_{\eta T}\right) }+ \\
		+\left\Vert \rho _{R,\min ,\lambda ,\xi }-\rho _{S}^{\ast }\right\Vert
		_{H^{2,1}\left( Q_{\eta T}\right) }\leq C_{1}\delta ^{1-\mu _{1}}, \\
		\left\Vert \beta _{\min ,\lambda ,\xi }-\beta ^{\ast }\right\Vert
		_{L_{2}\left( \Omega \right) }+\left\Vert \gamma _{\min ,\lambda ,\xi
		}-\gamma ^{\ast }\right\Vert _{L_{2}\left( \Omega \right) }\leq C_{1}\delta
		^{1-\mu _{1}},%
	\end{array}
	\right.  \label{5.27}
\end{equation}%
\emph{where the functions }$\rho _{S,\min ,\lambda ,\xi },\rho _{I,\min
	,\lambda ,\xi },\rho _{R,\min ,\lambda ,\xi }$\emph{\ are found using
	obvious analogs of (\ref{3.5}) and (\ref{3.10}), in which }$W$\emph{\ is
	replaced with }$\left( X_{\min ,\lambda ,\xi }+Y\right) \left( \mathbf{x}%
,t\right) $\emph{\ and coefficients }$\beta _{\min ,\lambda ,\xi }\left(
\mathbf{x}\right) $ \emph{and} $\gamma _{\min ,\lambda ,\xi }\left( \mathbf{x%
}\right) $\emph{\ are found via formulas (\ref{3.2}) and (\ref{3.3}), in
	which }$\rho _{S}$\emph{\ and }$\rho _{R}$\emph{\ are replaced with }$\rho
_{S,\min ,\lambda ,\xi }$\emph{\ and }$\rho _{R,\min ,\lambda ,\xi }$\emph{\
	respectively. }

\textbf{Theorem 5.3} (uniqueness). \emph{Suppose that conditions (\ref{3.1})
	hold. Then there exists at most one 5D vector function }$\Phi \in
C_{2}\left( \overline{\Omega }\right) \times C_{3}^{6,3}\left( \overline{Q}%
_{T}\right) $ \emph{in (\ref{2.11}) satisfying conditions (\ref{2.2})-(\ref%
	{2.5}), (\ref{2.9}) and (\ref{2.10}).}

\textbf{Proof of Theorem 5.3.} It is convenient to prove this theorem,
assuming that Theorem 5.2 is proven already. Assume that there exist two
vector functions $\Phi _{1},\Phi _{2}$ with the same data in (\ref{2.5}) and
(\ref{2.10}). Let $\overline{\Phi }=\Phi _{1}-\Phi _{2}.$ Setting $\delta =0$
in Theorem 5.2 and using (\ref{5.27}), we obtain $\overline{\Phi }\left(
\mathbf{x},t\right) =0$ in $Q_{\eta T}.$ To prove that $\overline{\Phi }%
\left( x,t\right) =0$ in $Q_{T},$ we refer to the well known uniqueness
result for parabolic equations with incomplete lateral Cauchy data \cite[%
Theorem 2.6.2]{KL}, \cite[Theorem 2 in \S 1 of Chapter 4 ]{LRS}. $\square $

\subsection{Global convergence of the gradient descent method}

\label{sec:5.4}

Assume that conditions of Theorem 5.2 hold. Suppose that
\begin{equation}
	W^{\ast }\in B^{\ast }\left( \frac{R}{3}\right) ,  \label{5.28}
\end{equation}%
where the set $B^{\ast }\left( R\right) $ is defined in (\ref{5.14}). Denote
\begin{equation}
	W_{\min ,\lambda ,\xi }=X_{\min ,\lambda ,\xi }+Y.  \label{5.280}
\end{equation}%
It follows from (\ref{5.2}), (\ref{5.27}) and (\ref{5.28}) that it is
reasonable to assume that
\begin{equation}
	W_{\min ,\lambda ,\xi }\in B\left( \frac{R}{3}\right) .  \label{5.29}
\end{equation}%
Consider an arbitrary point
\begin{equation}
	W_{0}\in B\left( \frac{R}{3}\right) .  \label{5.30}
\end{equation}%
Let $\sigma >0$ be a number. Consider the gradient descent method of the
minimization of the functional $J_{\lambda ,\xi }\left( W\right) ,$%
\begin{equation}
	W_{n}=W_{n-1}-\sigma J_{\lambda ,\xi }^{\prime }\left( W_{n-1}\right) ,\text{
	}n=1,2,...  \label{5.31}
\end{equation}%
Note that since by Theorem 5.1 $J_{\lambda ,\xi }^{\prime }\left( W\right)
\in H_{6,0}^{4}\left( Q_{T}\right) ,$ then (\ref{5.5}) and (\ref{5.31})
imply that all vector functions $W_{n}$ have the same boundary conditions as
ones in (\ref{3.11}).

\textbf{Theorem 5.4}. \emph{Assume that conditions of Theorem 5.2 as well as
	conditions (\ref{5.29})-(\ref{5.31}) hold. Then there exists a number }$%
\sigma _{0}\in \left( 0,1\right) $\emph{\ such that for any }$\sigma \in
\left( 0,\sigma _{0}\right) $\emph{\ there exists a number }$\theta \left(
\sigma \right) \in \left( 0,1\right) $\emph{\ such that for all }$n\geq 1$%
\begin{equation}
	W_{n}\in B\left( R\right)  \label{5.32}
\end{equation}%
and the following convergence estimates hold:%
\begin{equation*}
	\left\Vert W_{n}-W^{\ast }\right\Vert _{H^{2,1}\left( Q_{\gamma T}\right)
	}+\left\Vert \rho _{S,n}-\rho _{S}^{\ast }\right\Vert _{H^{2,1}\left(
		Q_{\gamma T}\right) }+\left\Vert \rho _{I,n}-\rho _{I}^{\ast }\right\Vert
	_{H^{2,1}\left( Q_{\gamma T}\right) }+
\end{equation*}%
\begin{equation}
	+\left\Vert \rho _{R,n}-\rho _{R}^{\ast }\right\Vert _{H^{2,1}\left(
		Q_{\gamma T}\right) }+\left\Vert \beta _{n}-\beta ^{\ast }\right\Vert
	_{L_{2}\left( \Omega \right) }+\left\Vert \gamma _{n}-\gamma ^{\ast
	}\right\Vert _{L_{2}\left( \Omega \right) }\leq  \label{5.33}
\end{equation}%
\begin{equation*}
	\leq C_{1}\delta ^{1-\mu _{1}}+C_{1}\theta ^{n}\left\Vert W_{\min ,\lambda
		,\xi }-W_{0}\right\Vert _{H_{6}^{4}\left( Q_{T}\right) },
\end{equation*}%
\emph{where functions }$\rho _{S,n},\rho _{I,n},\rho _{R,n}$\emph{\ and }$%
\beta _{n},\gamma _{n}$\emph{\ are obtained from }$W_{n} $\emph{\ by the
	obvious analogs of the procedures described in the formulation of Theorem
	5.2. }

\textbf{Proof of Theorem 5.4. }Suppose that Theorem 5.2 is proven. Relation (%
\ref{5.32}) as well as estimates%
\begin{equation}
	\left\Vert W_{\min ,\lambda ,\xi }-W_{n}\right\Vert _{H^{4}\left(
		Q_{T}\right) }\leq C_{1}\theta ^{n}\left\Vert W_{\min ,\lambda ,\xi
	}-W_{0}\right\Vert _{H^{4}\left( Q_{T}\right) },\text{ }n=1,2,...
	\label{5.34}
\end{equation}%
follow immediately from (\ref{5.29})-(\ref{5.31}) and \cite[Theorem 6]{SAR}.
Next, since (\ref{5.34}) holds, then (\ref{5.27}) and (\ref{5.280}) lead to
the target estimate (\ref{5.33}). $\square $

\textbf{Remark 5.1.}\emph{\ Since a smallness assumption is not imposed on
	the number }$R$\emph{\ and since the starting point }$W_{0}$ \emph{of
	iterations (\ref{5.31}) is an arbitrary point of the set }$B\left(
R/3\right) ,$\emph{\ then Theorem 5.4 establishes the global convergence of
	sequence (\ref{5.31}) to the true solution, as long as the noise level in
	the input data }$\delta \rightarrow 0,$\emph{\ see the last sentence of Step
	2 in section 1.}

\section{Proof of Theorem 5.1}

\label{sec:6}

Consider two arbitrary 6D vector functions $W_{1},W_{2}\in B\left( R\right) $
and denote
\begin{equation}
	h=W_{2}-W_{1}.  \label{1}
\end{equation}%
Then $W_{2}=W_{1}+h$. Also, (\ref{5.2}), (\ref{5.3}), (\ref{5.5}) and
triangle inequality imply:%
\begin{equation}
	\left.
	\begin{array}{c}
		h\in H_{6,0}^{4}\left( Q_{T}\right) \cap C_{6}^{2}\left( \overline{Q}
		_{T}\right) , \\
		\left\Vert h\right\Vert _{H_{6}^{4}\left( Q_{T}\right) }<2R,\text{ }
		\left\Vert h\right\Vert _{C_{6}^{2}\left( \overline{Q}_{T}\right) }\leq 2CR.%
	\end{array}
	\right.  \label{6.1}
\end{equation}%
Let $L\left( W\right) $ be the operator defined in (\ref{3.11}). Using the
multidimensional analog of Taylor's formula \cite[Theorem 4.8 of \S 4 of
Chapter 1]{V}, (\ref{1}) and (\ref{6.1}), we obtain
\begin{align}
	&L\left( W_{2}\right) =L\left( W_{1}+h\right) =L\left( W_{1}\right)
	+h_{t}-c\Delta h+A_{1,\text{lin}}\left( h,\mathbf{x,}t\right)  \notag \\
	& \hspace{1 cm} +A_{2,\text{lin}}\left( \nabla h,\mathbf{x,}t\right) +A_{3,%
		\text{lin}}\left( \int\limits_{T/2}^{t}h\left( \mathbf{x},\tau \right) d\tau
	,\mathbf{x,}t\right)  \label{6.2} \\
	& \hspace{0.5 cm} +A_{\text{nonlin}}\left( h,\nabla
	h,\int\limits_{T/2}^{t}h\left( \mathbf{x},\tau \right) d\tau ,\mathbf{x,}
	t\right) ,\text{ }\left( \mathbf{x,}t\right) \in Q_{T},  \notag
\end{align}
where the following estimates hold:%
\begin{equation*}
	\left\vert A_{1,\text{lin}}\left( h,\mathbf{x,}t\right) \right\vert
	,\left\vert A_{2,\text{lin}}\left( \nabla h,\mathbf{x,}t\right) \right\vert
	,\left\vert A_{3,\text{lin}}\left( \int\limits_{T/2}^{t}h\left( \mathbf{x}
	,\tau \right) d\tau ,\mathbf{x,}t\right) \right\vert \leq
\end{equation*}%
\begin{equation}
	\leq C_{1}\left( \left\vert h\right\vert \left( \mathbf{x,}t\right)
	+\left\vert \nabla h\right\vert \left( \mathbf{x,}t\right) +\left\vert
	\int\limits_{T/2}^{t}h\left( \mathbf{x},\tau \right) d\tau \right\vert
	\right) ,\text{ }\left( \mathbf{x,}t\right) \in Q_{T},  \label{6.3}
\end{equation}%
\begin{equation*}
	\left\vert A_{\text{nonlin}}\left( h,\nabla h,\int\limits_{T/2}^{t}h\left(
	\mathbf{x},\tau \right) d\tau ,\mathbf{x,}t\right) \right\vert \leq
\end{equation*}%
\begin{equation}
	\leq C_{1}\left[ h^{2}\left( \mathbf{x,}t\right) +\left\vert \nabla
	h\right\vert ^{2}\left( \mathbf{x,}t\right) +\left(
	\int\limits_{T/2}^{t}h\left( \mathbf{x},\tau \right) d\tau \right) ^{2} %
	\right] ,\text{ }\left( \mathbf{x,}t\right) \in Q_{T}.  \label{6.4}
\end{equation}%
In (\ref{6.2})-(\ref{6.4}) $A_{1,\text{lin}},A_{2,\text{lin}}$,$A_{3,\text{
		lin}}$ and $A_{\text{nonlin}}$ are linear and nonlinear operators
respectively with respect to $h$. Denote%
\begin{equation}
	A_{\text{lin}}=A_{1,\text{lin}}\left( h,\mathbf{x,}t\right) +A_{2,\text{lin}
	}\left( \nabla h,\mathbf{x,}t\right) +A_{3,\text{lin}}\left(
	\int\limits_{T/2}^{t}h\left( \mathbf{x},\tau \right) d\tau ,\mathbf{x,}
	t\right) ,  \label{6.5}
\end{equation}%
\begin{equation}
	A_{\text{nonlin}}=A_{\text{nonlin}}\left( h,\nabla
	h,\int\limits_{T/2}^{t}h\left( \mathbf{x},\tau \right) d\tau ,\mathbf{x,}
	t\right) .  \label{6.6}
\end{equation}%
Using (\ref{6.2})-(\ref{6.6}), we obtain%
\begin{align}
	&\left[ L\left( W_{1}+h\right) \right] ^{2}=\left[ L\left( W_{1}\right) %
	\right] ^{2}+2\left( h_{t}-c\Delta h+A_{\text{lin}}\right) L\left(
	W_{1}\right)  \notag \\
	& \hspace{0.5 cm} +2A_{\text{nonlin}}\cdot L\left( W_{1}\right) +\left(
	h_{t}-c\Delta h+A_{\text{lin}}+A_{\text{nonlin}}\right) ^{2}.  \label{6.7}
\end{align}

Let $\left[ ,\right] $ be the scalar product in $H_{6}^{4}\left(
Q_{T}\right) .$ It follows from (\ref{5.4}) and (\ref{6.7}) that%
\begin{equation*}
	J_{\lambda ,\xi }\left( W_{1}+h\right) -J_{\lambda ,\xi }\left( W_{1}\right)
	=
\end{equation*}%
\begin{equation*}
	=e^{-2\lambda b^{2}}\int\limits_{Q_{T}}2\left( h_{t}-c\Delta h+A_{\text{lin}
	}\right) L\left( W_{1}\right) \varphi _{\lambda }d\mathbf{x}dt+2\xi \left[
	W_{1},h\right] +
\end{equation*}%
\begin{align}
	&+e^{-2\lambda b^{2}}\int\limits_{Q_{T}}\left[ \left( h_{t}-c\Delta h+A_{
		\text{lin}}+A_{\text{nonlin}}\right) ^{2}+A_{\text{nonlin}}\cdot L\left(
	W_{1}\right) \right] \varphi _{\lambda }d\mathbf{x}dt  \label{6.8} \\
	&\hspace{2 cm} +\xi \left\Vert h\right\Vert _{H_{6}^{4}\left( Q_{T}\right)
	}^{2}.  \notag
\end{align}%
Consider the functional $\widetilde{J}_{\lambda ,\xi ,W_{1}}\left( h\right) $%
,%
\begin{equation}
	\widetilde{J}_{\lambda ,\xi ,W_{1}}\left( h\right) =e^{-2\lambda
		b^{2}}\int\limits_{Q_{T}}2\left( h_{t}-c\Delta h+A_{\text{lin}}\right)
	L\left( W_{1}\right) \varphi _{\lambda }d\mathbf{x}dt+2\xi \left[ W_{1},h %
	\right] .  \label{6.9}
\end{equation}%
By (\ref{6.1}) we can consider this functional as $\widetilde{J}_{\lambda
	,\xi }:H_{6,0}^{4}\left( Q_{T}\right) \rightarrow \mathbb{R}.$ Since $%
\widetilde{J}_{\lambda ,\xi }$ is a bounded linear functional, then Riesz
theorem implies that there exists a point $J_{\lambda ,\xi ,W_{1}}\in
H_{6,0}^{4}\left( Q_{T}\right) $ such that
\begin{equation}
	\widetilde{J}_{\lambda ,\xi }\left( h\right) =\left[ J_{\lambda ,\xi
		,W_{1}},h\right] ,\text{ }\forall h\in H_{6,0}^{4}\left( Q_{T}\right) .
	\label{6.10}
\end{equation}%
It follows from (\ref{6.3})-(\ref{6.10}) that
\begin{equation*}
	\lim_{\left\Vert h\right\Vert _{H_{6}^{4}\left( Q_{T}\right) }\rightarrow 0}
	\frac{J_{\lambda ,\xi }\left( W_{1}+h\right) -J_{\lambda ,\xi }\left(
		W_{1}\right) -\left[ J_{\lambda ,\xi ,W_{1}},h\right] }{\left\Vert
		h\right\Vert _{H_{6}^{4}\left( Q_{T}\right) }}=0.
\end{equation*}%
Hence, $J_{\lambda ,\xi ,W_{1}}$ is the Fr\'{e}chet derivative of the
functional $J_{\lambda ,\xi }\left( W\right) $ at the point $W_{1},$ i.e. $%
J_{\lambda ,\xi ,W_{1}}=J_{\lambda ,\xi }^{\prime }\left( W_{1}\right) \in
H_{6,0}^{4}\left( Q_{T}\right) $. Hence (\ref{6.8}) implies:%
\begin{equation*}
	J_{\lambda ,\xi }\left( W_{1}+h\right) -J_{\lambda ,\xi }\left( W_{1}\right)
	-J_{\lambda ,\xi }^{\prime }\left( W_{1}\right) \left( h\right) =
\end{equation*}%
\begin{equation}
	=e^{-2\lambda b^{2}}\int\limits_{Q_{T}}\left[ \left( h_{t}-c\Delta h+A_{
		\text{lin}}+A_{\text{nonlin}}\right) ^{2}+2A_{\text{nonlin}}\cdot L\left(
	W_{1}\right) \right] \varphi _{\lambda }d\mathbf{x}dt+  \label{6.11}
\end{equation}%
\begin{equation*}
	+\xi \left\Vert h\right\Vert _{H_{6}^{4}\left( Q_{T}\right) }^{2}.
\end{equation*}
We omit the proof of the Lipschitz stability properly (\ref{5.6}) of the
functional $J_{\lambda ,\xi }^{\prime }\left( W\right) $ since it is similar
with the proof of Theorem 3.1 of \cite{Bak}.

We now prove the key inequality (\ref{5.8}). Since $W_{1}\in B\left(
R\right) ,$ then the first line of (\ref{3.11}) and (\ref{5.1})-(\ref{5.3})
lead to $\left\vert L\left( W_{1}\right) \left( \mathbf{x},t\right)
\right\vert \leq C_{1}$ in $Q_{T}.$ Hence, using (\ref{6.3})-(\ref{6.6}) and
Cauchy-Schwarz inequality, we obtain the following estimate from the below
for the first term of the product in the integrand in the second line of (%
\ref{6.11}):%
\begin{equation*}
	\left( h_{t}-c\Delta h+A_{\text{lin}}+A_{\text{nonlin}}\right) ^{2}\left(
	\mathbf{x,}t\right) +2A_{\text{nonlin}}\left( \mathbf{x,}t\right) \cdot
	L\left( W_{1}\right) \left( \mathbf{x,}t\right) \geq
\end{equation*}%
\begin{equation*}
	\geq \frac{1}{2}\left( h_{t}-c\Delta h\right) ^{2}\left( \mathbf{x,}t\right)
	-C_{1}\left[ h^{2}\left( \mathbf{x,}t\right) +\left\vert \nabla h\right\vert
	^{2}\left( \mathbf{x,}t\right) +\left( \int\limits_{T/2}^{t}h\left( \mathbf{%
		x },\tau \right) d\tau \right) ^{2}\right] .
\end{equation*}%
Substitute this in (\ref{6.11}) and use the Carleman estimate of Theorem
4.1. We obtain%
\begin{equation*}
	J_{\lambda ,\xi }\left( W_{1}+h\right) -J_{\lambda ,\xi }\left( W_{1}\right)
	-J_{\lambda ,\xi }^{\prime }\left( W_{1}\right) \left( h\right) \geq
\end{equation*}%
\begin{equation*}
	\geq \frac{C}{\lambda }e^{-2\lambda
		b^{2}}\int\limits_{Q_{T}}h_{t}^{2}\varphi _{\lambda }d\mathbf{x}dt+\frac{C}{
		\lambda }e^{-2\lambda
		b^{2}}\sum\limits_{i,j=1}^{2}\int\limits_{Q_{T}}h_{x_{i}x_{j}}^{2}\varphi
	_{\lambda }d\mathbf{x}dt+
\end{equation*}%
\begin{equation*}
	+Ce^{-2\lambda b^{2}}\int\limits_{Q_{T}}\left( \lambda \left( \nabla
	h\right) ^{2}+\lambda ^{3}h^{2}\right) \varphi _{\lambda }d\mathbf{x}dt-
\end{equation*}%
\begin{equation*}
	-C_{1}e^{-2\lambda b^{2}}\int\limits_{Q_{T}}\left[ h^{2}\left( \mathbf{x,}
	t\right) +\left\vert \nabla h\right\vert ^{2}\left( \mathbf{x,}t\right)
	+\left( \int\limits_{T/2}^{t}h\left( \mathbf{x},\tau \right) d\tau \right)
	^{2}\right] \varphi _{\lambda }d\mathbf{x}dt+\xi \left\Vert h\right\Vert
	_{H_{6}^{4}\left( Q_{T}\right) }^{2}-
\end{equation*}%
\begin{equation*}
	-C\left( \left\Vert h\left( \mathbf{x},T\right) \right\Vert _{H^{1}\left(
		\Omega \right) }^{2}+\left\Vert h\left( \mathbf{x},0\right) \right\Vert
	_{H^{1}\left( \Omega \right) }^{2}\right) \lambda ^{2}\exp \left( -\lambda
	\frac{T^{2}}{2}\right) ,\text{ }\lambda \geq \lambda _{0}.
\end{equation*}%
Hence, using Lemma 4.1, we obtain that there exists a sufficiently large
number $\lambda _{1}=\lambda _{1}\left( \kappa ,c,R,Q_{T}\right) \geq
\lambda _{0}$ such that%
\begin{equation*}
	J_{\lambda ,\xi }\left( W_{1}+h\right) -J_{\lambda ,\xi }\left( W_{1}\right)
	-J_{\lambda ,\xi }^{\prime }\left( W_{1}\right) \left( h\right) \geq
\end{equation*}%
\begin{equation}
	\geq \frac{C_{1}}{\lambda }e^{-2\lambda
		b^{2}}\int\limits_{Q_{T}}h_{t}^{2}\varphi _{\lambda }d\mathbf{x}dt+\frac{
		C_{1}}{\lambda }e^{-2\lambda
		b^{2}}\sum\limits_{i,j=1}^{2}\int\limits_{Q_{T}}h_{x_{i}x_{j}}^{2}\varphi
	_{\lambda }d\mathbf{x}dt+  \label{6.12}
\end{equation}%
\begin{equation*}
	+C_{1}e^{-2\lambda b^{2}}\int\limits_{Q_{T}}\left( \lambda \left( \nabla
	h\right) ^{2}+\lambda ^{3}h^{2}\right) \varphi _{\lambda }d\mathbf{x}dt+
\end{equation*}%
\begin{equation*}
	+\xi \left\Vert h\right\Vert _{H_{6}^{4}\left( Q_{T}\right)
	}^{2}-C_{1}\left( \left\Vert h\left( \mathbf{x},T\right) \right\Vert
	_{H^{1}\left( \Omega \right) }^{2}+\left\Vert h\left( \mathbf{x},0\right)
	\right\Vert _{H^{1}\left( \Omega \right) }^{2}\right) \lambda ^{2}\exp
	\left( -\lambda \frac{T^{2}}{2}\right) ,\text{ }\forall \lambda \geq \lambda
	_{1}.
\end{equation*}%
By trace theorem and (\ref{5.7})%
\begin{equation}
	\frac{\xi }{2}\left\Vert h\right\Vert _{H_{6}^{4}\left( Q_{T}\right)
	}^{2}-C_{1}\left( \left\Vert h\left( \mathbf{x},T\right) \right\Vert
	_{H^{1}\left( \Omega \right) }^{2}+\left\Vert h\left( \mathbf{x},0\right)
	\right\Vert _{H^{1}\left( \Omega \right) }^{2}\right) \lambda ^{2}\exp
	\left( -\lambda \frac{T^{2}}{2}\right) \geq 0,  \label{6.13}
\end{equation}%
for all $\lambda \geq \lambda _{1}.$ The target estimate (\ref{5.8}) follows
immediately from (\ref{4.20}), (\ref{1}), (\ref{6.12}) and (\ref{6.13}).

The existence and uniqueness of the minimizer $W_{\min ,\lambda ,\xi }$ of
the functional $J_{\lambda ,\xi }\left( W\right) $ on the set $\overline{
	B\left( R\right) }$ as well as inequality (\ref{5.9}) follow from a
combination of Lemma 2.1 with Theorem 2.1 of \cite{Bak} as well as from a
combination of Lemma 5.2.1 with Theorem 5.2.1 of \cite{KL}. \ $\square $

\section{Proof of Theorem 5.2}

\label{sec:7}

Recall that we work in this theorem with $\lambda \geq \lambda _{2},$ where $%
\lambda _{2}$ is defined in (\ref{5.22}), and the regularization parameter $%
\xi $ is chosen now as in (\ref{5.23}). Also, recall that the direct analog
of Theorem 5.1 is valid for the functional $I_{\lambda ,\xi }\left( X\right)
=J_{\lambda ,\xi }\left( X+Y\right) $ on the set $\overline{B_{0}\left(
	2R\right) }.$ This functional is defined in (\ref{5.21}), and the set $%
B_{0}\left( R\right) $ is defined in (\ref{5.18}). Finally, we recall that $%
X_{\min ,\lambda ,\xi }\in \overline{B_{0}\left( 2R\right) }$ is the unique
minimizer of the functional $I_{\lambda ,\xi }\left( X\right) $ on the set $%
\overline{B_{0}\left( 2R\right) }.$

Using (\ref{5.19}) and the obvious analog of estimate (\ref{5.8}), we obtain%
\begin{equation}
	\left.
	\begin{array}{c}
		I_{\lambda ,\xi }\left( \widetilde{W^{\ast }}\right) -I_{\lambda ,\xi
		}\left( X_{\min ,\lambda ,\xi }\right) -I_{\lambda ,\xi }\left( X_{\min
			,\lambda ,\xi }\right) \left( \widetilde{W^{\ast }}-X_{\min ,\lambda ,\xi
		}\right) \geq \\
		\geq \left( C_{1}/\lambda \right) \exp \left( -\lambda \eta
		^{2}T^{2}/2\right) \left\Vert \widetilde{W^{\ast }}-X_{\min ,\lambda ,\xi
		}\right\Vert _{H^{2,1}\left( Q_{\eta T}\right) }^{2},\text{ }\forall \lambda
		\geq \lambda _{2}.%
	\end{array}
	\right.  \label{7.1}
\end{equation}%
Since $I_{\lambda ,\xi }\left( X_{\min ,\lambda ,\xi }\right) \geq 0,$ then (%
\ref{5.9}) implies:
\begin{equation*}
	-I_{\lambda ,\xi }\left( X_{\min ,\lambda ,\xi }\right) -I_{\lambda ,\xi
	}\left( X_{\min ,\lambda ,\xi }\right) \left( \widetilde{W^{\ast }}-X_{\min
		,\lambda ,\xi }\right) \leq 0.
\end{equation*}%
Hence, by (\ref{7.1})%
\begin{equation}
	\left\Vert \widetilde{W^{\ast }}-X_{\min ,\lambda ,\xi }\right\Vert
	_{H^{2,1}\left( Q_{\eta T}\right) }^{2}\leq C_{1}\lambda \exp \left( \lambda
	\frac{\eta ^{2}T^{2}}{2}\right) I_{\lambda ,\xi }\left( \widetilde{W^{\ast }}
	\right) .  \label{7.2}
\end{equation}%
Using (\ref{5.4}), (\ref{5.19}), (\ref{5.21}) and (\ref{5.23}), we obtain
\begin{equation}
	\left.
	\begin{array}{c}
		I_{\lambda ,\xi }\left( \widetilde{W^{\ast }}\right) =e^{-2\lambda
			b^{2}}\int\limits_{Q_{T}}\left( L\left( W^{\ast }-Y^{\ast }+Y\right) \right)
		^{2}\varphi _{\lambda }d\mathbf{x}dt+ \\
		+2\exp \left( -\lambda T^{2}/4\right) \left\Vert W^{\ast }-Y^{\ast
		}+Y\right\Vert _{H_{6}^{4}\left( Q_{T}\right) }^{2}.%
	\end{array}
	\right.  \label{7.3}
\end{equation}%
By (\ref{5.13})
\begin{equation}
	\widetilde{L}\left( W^{\ast }\right) =0.  \label{7.4}
\end{equation}%
The operator $\widetilde{L}$ is obtained from the operator $L$ via replacing
in (\ref{3.6}) the vector function $\left( p_{1},p_{2},p_{3}\right) \left(
\mathbf{x}\right) $ with the vector function $\left( p_{1}^{\ast
},p_{2}^{\ast },p_{3}^{\ast }\right) \left( \mathbf{x}\right) .$ Hence,
using (\ref{5.15}), (\ref{5.17}) and (\ref{7.4}), we obtain
\begin{equation*}
	\left.
	\begin{array}{c}
		\left( L\left( W^{\ast }-Y^{\ast }+Y\right) \right) ^{2}=\left[ \widetilde{L}
		\left( W^{\ast }\right) +\left( L\left( W^{\ast }-Y^{\ast }+Y\right) -
		\widetilde{L}\left( W^{\ast }\right) \right) \right] ^{2}= \\
		=\left( L\left( W^{\ast }-Y^{\ast }+Y\right) -\widetilde{L}\left( W^{\ast
		}\right) \right) ^{2}\leq C_{1}\delta ^{2}.%
	\end{array}
	\right.
\end{equation*}%
Hence, by (\ref{4.2}), (\ref{5.23}) and (\ref{7.3})%
\begin{equation*}
	I_{\lambda ,\xi }\left( \widetilde{W^{\ast }}\right) \leq C_{1}\delta
	^{2}+C_{1}\exp \left( -\lambda \frac{T^{2}}{4}\right) .
\end{equation*}%
Hence, (\ref{7.2}) implies%
\begin{equation}
	\left\Vert \widetilde{W^{\ast }}-X_{\min ,\lambda ,\xi }\right\Vert
	_{H^{2,1}\left( Q_{\eta T}\right) }\leq C_{1}\left( \sqrt{\lambda }\delta
	+\exp \left( -\lambda T^{2}/8\right) \right) \exp \left( \lambda \frac{\eta
		^{2}T^{2}}{4}\right) .  \label{7.5}
\end{equation}%
Choose $\lambda =\lambda \left( \delta \right) $ such that $\exp \left(
-\lambda \left( \delta \right) T^{2}/8\right) =\delta .$ Hence, (\ref{5.24})
should be in place. Hence, using (\ref{7.5}), we obtain%
\begin{equation}
	\left\Vert \widetilde{W^{\ast }}-X_{\min ,\lambda ,\xi }\right\Vert
	_{H^{2,1}\left( Q_{\eta T}\right) }\leq C_{1}\sqrt{\lambda \left( \delta
		\right) }\exp \left( \lambda \left( \delta \right) \frac{\eta ^{2}T^{2}}{4}
	\right) \delta =C_{1}\sqrt{\lambda \left( \delta \right) }\delta ^{1-2\eta
		^{2}}.  \label{7.6}
\end{equation}%
Choose $\eta ^{2}$ as in (\ref{5.25}). Then in (\ref{7.6}) $\sqrt{\lambda
	\left( \delta \right) }\delta ^{1-2\eta ^{2}}=\sqrt{\lambda \left( \delta
	\right) }\delta ^{1-\mu _{1}+\mu _{2}}.$ We can choose $\delta _{0}$ in (\ref%
{5.24}) so small that $\sqrt{\lambda \left( \delta \right) }\delta ^{1-\mu
	_{1}+\mu _{2}}<\delta ^{1-\mu _{1}}$ for all $\delta \in \left( 0,\delta
_{0}\right) .$

Hence, (\ref{5.17}), (\ref{5.19}), (\ref{7.6}) and triangle inequality imply%
\begin{equation*}
	\left.
	\begin{array}{c}
		C_{1}\delta ^{1-\mu _{1}}\geq \left\Vert \widetilde{W^{\ast }}-X_{\min
			,\lambda ,\xi }\right\Vert _{H^{2,1}\left( Q_{\gamma T}\right) } \\
		=\left\Vert W^{\ast }-\left( X_{\min ,\lambda ,\xi }+Y\right) -\left(
		Y^{\ast }-Y\right) \right\Vert _{H^{2,1}\left( Q_{\gamma T}\right) }\geq \\
		\geq \left\Vert \left( X_{\min ,\lambda ,\xi }+Y\right) -W^{\ast
		}\right\Vert _{H^{2,1}\left( Q_{\gamma T}\right) }-\delta .%
	\end{array}
	\right.
\end{equation*}%
Hence,
\begin{equation}
	\left\Vert \left( X_{\min ,\lambda ,\xi }+Y\right) -W^{\ast }\right\Vert
	_{H^{2,1}\left( Q_{\gamma T}\right) }\leq C_{1}\delta ^{1-\mu _{1}}+\delta .
	\label{7.7}
\end{equation}%
Since $\delta <\delta ^{1-\mu _{1}},$ then (\ref{7.7}) implies the estimate
in the first line of (\ref{5.27}). The other two estimates of (\ref{5.27})
follow. $\ \ \square $

\section{Numerical Studies}

\label{sec:8}

In this section we describe our numerical studies. We specify the domain $%
\Omega $ in \eqref{2.1} and the number $T$ as
\begin{equation}
	\Omega =\left\{ \mathbf{x}=\left( x,y\right) :a<x<b,\left\vert y\right\vert
	<A\right\} ,\ a=0.1,b=1.1,A=0.5,\ T=1.  \label{8.1}
\end{equation}%
We took $\eta _{S}=\eta _{I}=\eta _{R}=0.01$ in \eqref{2.7}. The velocities
of movements of the populations were chosen as:
\begin{equation}
	q_{S}\left( \mathbf{x}\right) =q_{I}\left( \mathbf{x}\right) =q_{R}\left(
	\mathbf{x}\right) =
	\begin{pmatrix}
		0.2 \\
		0.2%
	\end{pmatrix}
	.  \label{8.2}
\end{equation}%
We have assigned zero fluxes of the populations through the boundary in %
\eqref{2.5},
\begin{equation}
	g_{1}(\mathbf{x},t)=g_{2}(\mathbf{x},t)=g_{3}(\mathbf{x},t)=0.  \label{8.3}
\end{equation}%
The initial conditions in \eqref{2.6} were chosen as:
\begin{align}
	& \rho _{S}^{0}(\mathbf{x})=0.6\exp \left( -10\left(
	(x-0.6)^{2}+y^{2}\right) \right) ,  \label{8.4} \\
	& \rho _{I}^{0}(\mathbf{x})=0.6\exp \left( -35\left(
	(x-0.7)^{2}+(y-0.1)^{2}\right) \right) ,  \label{8.5} \\
	& \rho _{R}^{0}(\mathbf{x})=0.  \label{8.6}
\end{align}

\textbf{Remark 8.1.} \emph{To better demonstrate the effectiveness of our
	numerical method, we choose letter-like shapes of abnormalities to be
	imaged. Indeed, letters have non-convex shapes with voids. We want to rather
	accurately image both shapes of abnormalities and values of unknown
	coefficients }$\beta (\mathbf{x})$\emph{\ and }$\gamma (\mathbf{x})$ \emph{\
	inside and outside of them}$.$\emph{\ However, we are not particularly
	concerned with accurate imaging of sharp edges of abnormalities. }

The infection rate $\beta (\mathbf{x})$ was chosen as:
\begin{equation}
	\beta (\mathbf{x})=\left\{
	\begin{array}{cc}
		c_{\beta }=const.>0, & \mbox{inside the tested inclusion,} \\
		0.1, & \mbox{outside the tested inclusion.}%
	\end{array}
	\right.  \label{8.7}
\end{equation}%
In the numerical tests below, we take $c_{\beta }=0.6,0.4$, and the
inclusions with the shapes of the letters `$A$' and `$B$'. The recovery rate
was chosen as:
\begin{equation}
	\gamma (\mathbf{x})=\left\{
	\begin{array}{cc}
		c_{\gamma }=const.>0, & \mbox{inside the tested inclusion,} \\
		0.1, & \mbox{outside the tested inclusion.}%
	\end{array}
	\right.  \label{8.8}
\end{equation}%
In the numerical tests below, we take $c_{\gamma }=0.4,0.6$, and the
inclusions with the shapes of the letters `$\Omega $' and `$D$'. Hence,
inclusion/background contrasts for both unknown coefficients are 6:1 and 4:1.

We used the classic implicit scheme with the spatial mesh sizes $1/80\times
1/80$ and the temporal mesh step size $1/320$ to solve the forward problem %
\eqref{2.2}-\eqref{2.6} for data generation. In the computations of the
Minimization Problem, the spatial mesh sizes were $1/20\times 1/20$ and the
temporal mesh step size was $1/10$. The integral and differential operators
in the Minimization Problem \eqref{5.4} are been approximated by the
trapezoidal formula and the finite differences respectively with respect to
the values of 6D vector function $W$ at the grid points. The minimization
was performed with respect to the values of this vector function at grid
points.

To guarantee that the solution of the problem of the minimization of the
functional $J_{\lambda ,\xi }(W)$ in \eqref{5.4} satisfies the boundary
conditions in \eqref{3.11}, we adopt the Matlab's built-in optimization
toolbox \textbf{fmincon} to minimize the discretized form of the functional $%
J_{\lambda ,\xi }(W)$. The regularization parameter $\xi =0.01$ in %
\eqref{5.4} was chosen by the trial and error procedure. The iterations of
\textbf{fmincon} stop when the following condition is met:
\begin{equation}
	|\nabla J_{\lambda ,\xi }(W)|<10^{-2}.  \label{8.9}
\end{equation}%
The Neumann boundary conditions in \eqref{3.11} were approximated by the
finite difference method on the boundary in the iterations of \textbf{fmincon%
}. Denote the discrete points
\begin{equation}
	\begin{split}
		& \mathbf{x}_{i,j}=(x_{i},y_{j})=(a+ih_{x},-A+jh_{y}),\ i=0,1,\cdots
		,N_{x},\ j=0,1,\cdots ,N_{y}, \\
		& \hspace{2cm}N_{x}=(b-a)/h_{x},\ N_{y}=2A/h_{y},\ h_{x}=h_{y}=1/20.
	\end{split}
	\label{8.901}
\end{equation}%
Then the iterations of \textbf{fmincon} have the following constraint
conditions:
\begin{align}
	& 4W(x,y_{1},t)-W(x,y_{0},t)-3W(x,y_{2},t)=2h_{y}G(x,-A,t),  \label{8.902} \\
	& -4W(x,y_{N_{y}-1},t)+W(x,y_{N_{y}-2},t)+3W(x,y_{N_{y}},t)=2h_{y}G(x,A,t),
	\label{8.903} \\
	& 4W(x_{1},y,t)-W(x_{0},y,t)-3W(x_{2},y,t)=2h_{x}G(a,y,t),  \label{8.904} \\
	& -4W(x_{N_{x}-1},y,t)+W(x_{N_{x}-2},y,t)=2h_{x}G(b,y,t)-3F(b,y,t).
	\label{8.905}
\end{align}%
The formula \eqref{8.905} also contains the Dirichlet boundary conditions in %
\eqref{3.11} as
\begin{equation}
	W(\mathbf{x},t)\mid _{\Gamma _{T}}=F(\mathbf{x},t).  \label{8.906}
\end{equation}

Using the Dirichlet boundary conditions in \eqref{3.11}, we use the linear
interpolation inside the domain $Q_{T}$ in the $x-$direction to calculate
the starting point $W^{0}$ of our iterations as:
\begin{equation}
	W^{0}(\mathbf{x},t)=\frac{x-a}{b-a}F(\mathbf{x},t).  \label{8.907}
\end{equation}
In other words, we use linear interpolation inside the domain $Q_{T}$ of the
Dirichlet boundary conditions \eqref{3.11}. Although it follows from (\ref%
{8.907}) that the starting point $W^{0}(\mathbf{x},t)$ satisfies only
Dirichlet boundary conditions in \eqref{3.11} and does not satisfy Neumann
boundary conditions in \eqref{3.11}, still \eqref{8.902}-\eqref{8.906} imply
that boundary conditions in \eqref{3.11} are satisfied on all other
iterations of \textbf{fmincon}.

When we consider the noisy pollution in the observation data in \eqref{2.9}
and \eqref{2.10}, the random noisy is given as
\begin{equation*}
	p_{i,\zeta }(\mathbf{x})=p_{i}(\mathbf{x})(1+\zeta _{p_{i}}),\ f_{i,\zeta
	}(y,t)=f_{i}(y,t)(1+\sigma \zeta _{f_{i},y,t}),\ i=1,2,3,
\end{equation*}%
where $\zeta _{p_{i}},i=1,2,3$ are the uniformly distributed random
variables in the interval $[-1,1]$ depending on the point $\mathbf{x}\in
\Omega $, $\zeta _{f_{i},y,t},i=1,2,3$ are the uniformly distributed random
variables in the interval $[-1,1]$ depending on the point $y\in \left[
-0.5,0.5\right] $ and time $t\in \lbrack 0,T]$, and $\sigma
=0.01,0.02,0.03,0.05$ are the noisy level, which correspond to $1\%,2\%,3\%$
and $5\%$ noise levels respectively. To numerically calculate the first and
second derivatives of noisy functions $p_{i,\zeta }$ and $f_{i,\zeta }$, $%
i=1,2,3$, we use the natural cubic splines to approximate the noisy
observation data, and then use the derivatives of those splines to
approximate the derivatives of corresponding noisy observation data. With
the spatial mesh sizes $1/20\times 1/20$ and the temporal mesh step size $%
1/10$, we generate the cubic splines $s_{p,i},s_{f,i}$ of functions $%
p_{i,\zeta }$, $f_{i,\zeta }$, and then calculate the first and second
derivatives of $s_{p,i},s_{f,i}$ to approximate the first and second
derivatives respect to $x,y,t$ of $p_{i,\zeta }$, $f_{i,\zeta }$, $i=1,2,3$.

\begin{figure}[tbp]
	\includegraphics[width = 4.8in]{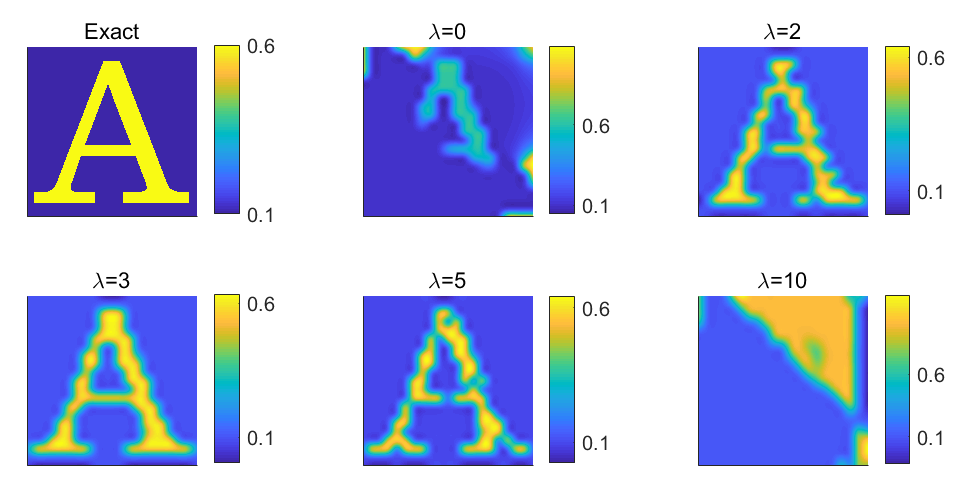}
	\caption{Test 1. The reconstructed results of function $\protect\beta(
		\mathbf{x})$ with different $\protect\lambda$ in \eqref{5.4}. The function $%
		\protect\beta(\mathbf{x})$ is given in \eqref{8.7} with $c_{\protect\beta %
		}=0.6$ inside of the letter `A' and the function $\protect\gamma(\mathbf{x})$
		is given in \eqref{8.8} with $c_{\protect\gamma}=0.4$ inside of the letter `$%
		\Omega$'.}
	\label{plot_diff_lambda_beta}
\end{figure}

\indent\textbf{Test 1.} We test the case when the inclusion of $\beta (%
\mathbf{x})$ in \eqref{8.7} has the shape of the letter `$A$' with $c_{\beta
}=0.6$ and the inclusion of $\gamma (\mathbf{x})$ in \eqref{8.8} has the
shape of the letter `$\Omega $' with $c_{\gamma }=0.4$. First, we work on
the reference case to figure out an optimal value of the parameter $\lambda $
of the Carleman Weight Function $\varphi _{\lambda }\left( \mathbf{x}%
,t\right) $ in (\ref{4.1}). The results with different values of $\lambda $
are displayed on Figure \ref{plot_diff_lambda_beta} and Figure \ref%
{plot_diff_lambda_gamma}. One can observe that the images have a low quality
with too small or too large $\lambda $. The low quality at $\lambda =0$
means that the \ CWF should be present indeed in the functional $J_{\lambda
	,\xi }\left( W\right) .$ Thus, we choose $\lambda =3$ as the optimal value
of this parameter.

\begin{figure}[htbp]
	\centering
	\includegraphics[width = 5in]{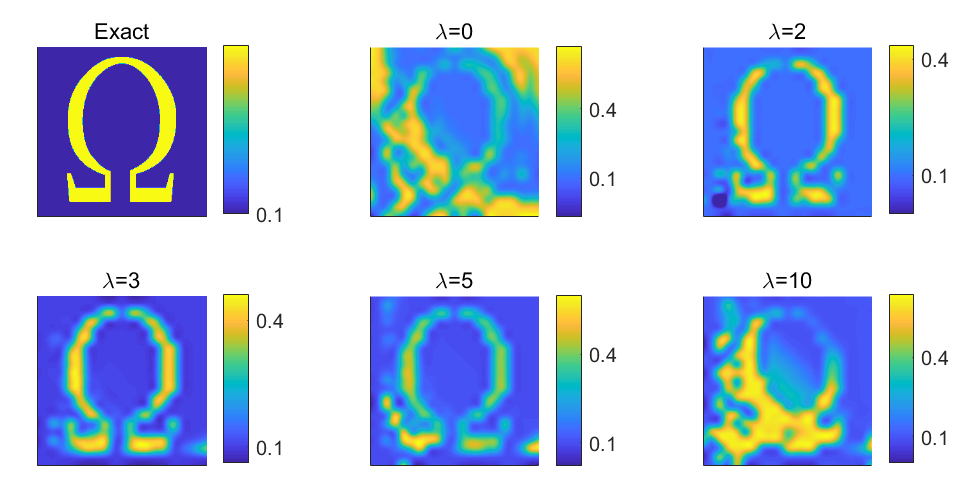}
	\caption{Test 1. The reconstructed results of function $\protect\gamma(
		\mathbf{x})$ with different $\protect\lambda$ in \eqref{5.4}. The function $%
		\protect\beta(\mathbf{x})$ is given in \eqref{8.7} with $c_{\protect\beta %
		}=0.6$ inside of the letter `A' and the function $\protect\gamma(\mathbf{x})$
		is given in \eqref{8.8} with $c_{\protect\gamma}=0.4$ inside of the letter `$%
		\Omega$'. Thus, it is clear from Figure \protect\ref{plot_diff_lambda_beta}
		and Figure \protect\ref{plot_diff_lambda_gamma} that the presence of the CWF
		in the functional $J_{\protect\lambda ,\protect\xi }\left( W\right) $ is
		necessary since the images for $\protect\lambda =0$ are poor. It is also
		clear that $\protect\lambda =3$ is the optimal value of the parameter $%
		\protect\lambda$.}
	\label{plot_diff_lambda_gamma}
\end{figure}

\begin{figure}[htbp]
	\centering
	\includegraphics[width = 5in]{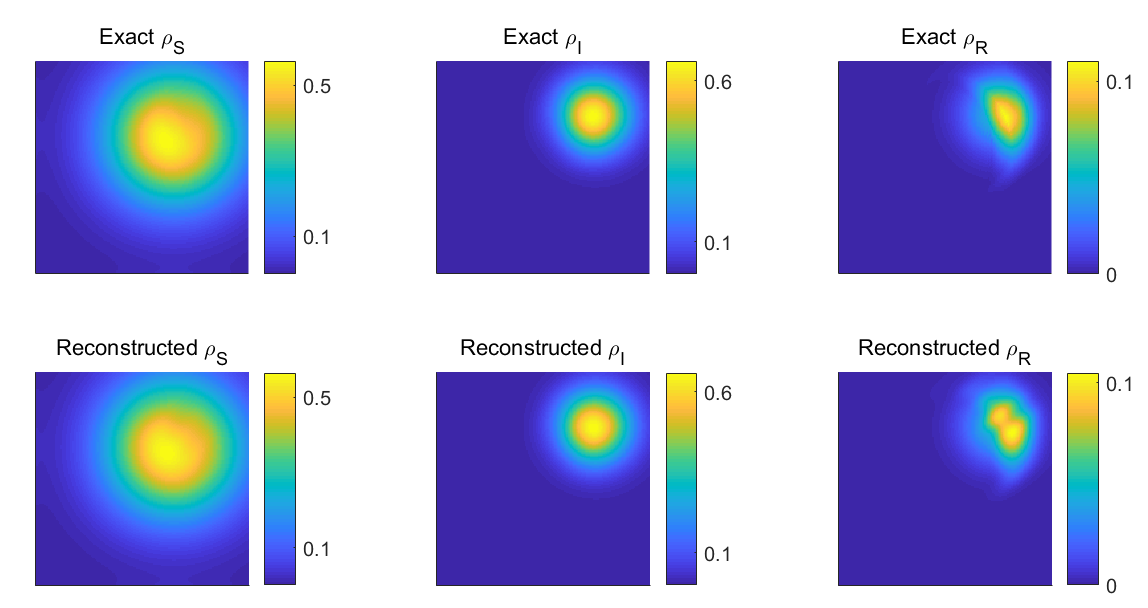}
	\caption{Test 1. The exact (top line) and reconstructed (bottom line)
		functions $\protect\rho_{S}$ (left column), $\protect\rho_{I}$ (middle line)
		and $\protect\rho_{R}$ (right line) at $t=0.75$, where the function $\protect%
		\beta(\mathbf{x})$ is given in \eqref{8.7} with $c_{\protect\beta}=0.6$
		inside of the letter `A' and the function $\protect\gamma(\mathbf{x})$ is
		given in \eqref{8.8} with $c_{\protect\gamma}=0.4$ inside of the letter `$%
		\Omega$'. }
	\label{plot_re_SIR}
\end{figure}

Using the solution of the Minimization Problem, we can reconstruct the
functions $\rho _{S},\rho _{I},\rho _{R}$ via \eqref{3.5}. At $t=0.75$, the
solutions of the forward problem \eqref{2.2}-\eqref{2.6} with known
coefficient functions $\beta (\mathbf{x})$ and $\gamma (\mathbf{x})$ are
indicated by `Exact', and the reconstructed functions $\rho _{S},\rho
_{I},\rho _{R}$ by \eqref{3.5} are indicated by `Reconstructed' in Figure %
\ref{plot_re_SIR}. The reconstructed functions $\rho _{S},\rho _{I},\rho
_{R} $ are almost as same as the solutions of the forward problem. In other
words, solving our Coefficient Inverse Problem \eqref{2.9}, \eqref{2.10} we
can accurately calculate the values of functions $\rho _{S}\left( \mathbf{x}%
,t\right) ,\rho _{I}\left( \mathbf{x},t\right) $ and $\rho _{R}\left(
\mathbf{x},t\right) $ for $\left( \mathbf{x},t\right) \in Q_{T}.$

\textbf{Test 2.} We test the case when the inclusion of $\beta (\mathbf{x})$
in \eqref{8.7} has the shape of the letter `$B$' with $c_{\beta }=0.4$ and
the inclusion of $\gamma (\mathbf{x})$ in \eqref{8.8} has the shape of the
letter `$D$' with $c_{\gamma }=0.6$. The result is displayed in Figure \ref%
{plot_re_beta_B_gamma_D}.
\begin{figure}[htbp]
	\centering
	\includegraphics[width = 3in]{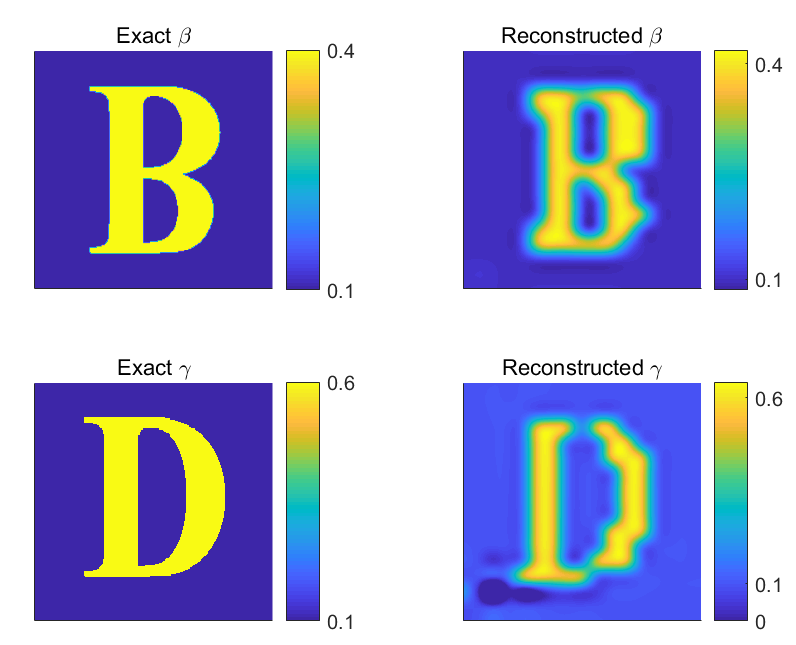}
	\caption{Test 2. The exact (left) and reconstructed (right) functions $%
		\protect\beta (\mathbf{x})$ (top line) and $\protect\gamma (\mathbf{x})$
		(bottom line), where the function $\protect\beta (\mathbf{x})$ is given in
		\eqref{8.7} with $c_{\protect\beta }=0.4$ inside of the letter `B' and the
		function $\protect\gamma (\mathbf{x})$ is given in \eqref{8.8} with $c_{
			\protect\gamma }=0.6$ inside of the letter `D'. }
	\label{plot_re_beta_B_gamma_D}
\end{figure}

\begin{figure}[htbp]
	\centering
	\includegraphics[width = 5in]{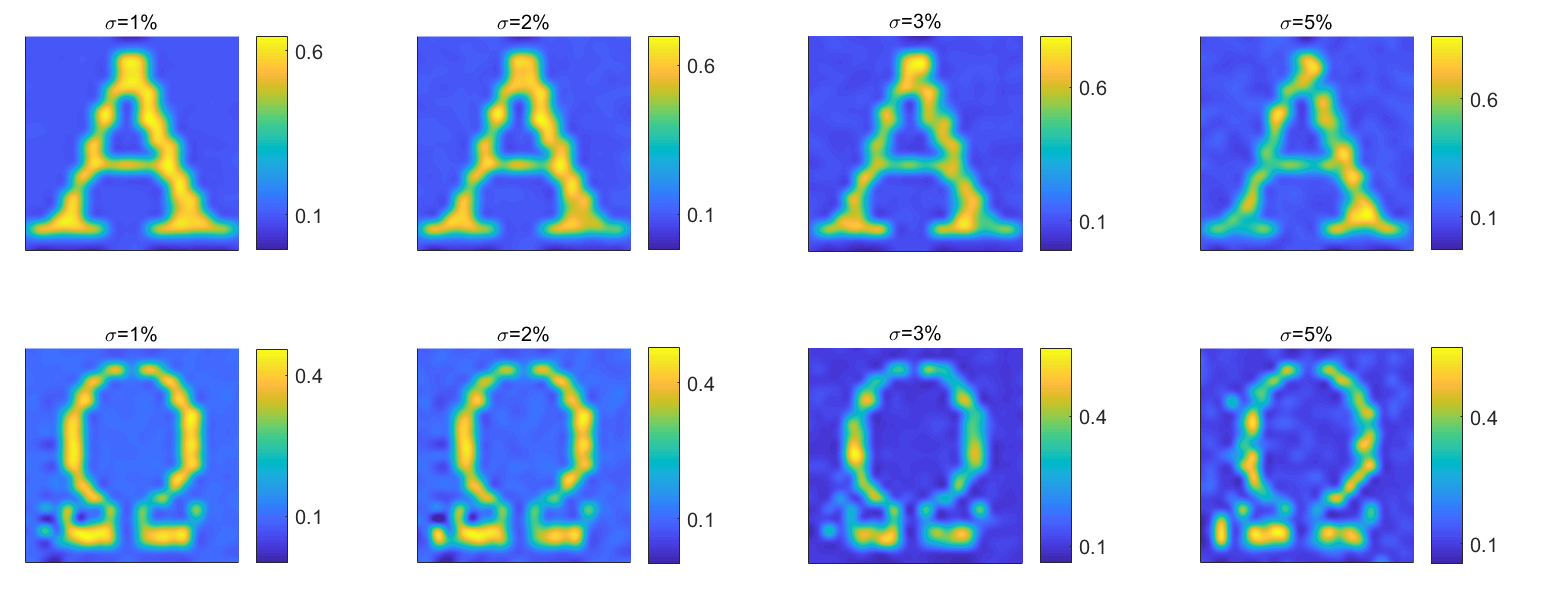}
	\caption{Test 3. The reconstructed functions $\protect\beta (\mathbf{x})$
		(top line) and $\protect\gamma (\mathbf{x})$ (bottom line) with $\protect%
		\sigma =0.01$ (1st column), $\protect\sigma =0.02$ (2ed column), $\protect%
		\sigma =0.03$ (3rd column) and $\protect\sigma =0.05$ (4th column). }
	\label{plot_beta_A_gamma_Omega_noisy}
\end{figure}

\textbf{Test 3.} Following \textbf{Test 1}, the noisy data is considered.
The result is displayed in Figure \ref{plot_beta_A_gamma_Omega_noisy}.

\section{Summary}

\label{sec:9}

This is the first work, \ in which a Coefficient Inverse Problem for the
spatiotemporal SIR model of \cite{Lee} is considered. We have constructed a
new version of the convexification method for this CIP. The main conclusion
of our convergence analysis is the global convergence of this method. In
other words, its convergence to the true solution is guaranteed regardless
on the availability of a good initial guess. The key to the global
convergence property is the presence of the Carleman Weight Function in the
resulting Tikhonov-like regularization functional, and this is clearly
confirmed in the numerical Test 1, see Figure \ref{plot_diff_lambda_beta}
and Figure \ref{plot_diff_lambda_gamma}. As a by-product, uniqueness theorem
for our CIP is proven.

Results of our numerical studies demonstrate an accurate performance of our
method for rather complicated shapes of abnormalities. Since the minimal
amount of measured data are required in our CIP, then results of this paper
indicate a possibility of a significant decrease of the cost of the
spatiotemporal monitoring of the spread of epidemics in affected cities.


\end{document}